\documentclass[12pt,a4paper]{amsart}
\usepackage{pinlabel}
\usepackage{amssymb}

\input xy
\xyoption{all}

\numberwithin{equation}{section}
\numberwithin{equation}{subsection}

\theoremstyle{plain}
\newtheorem{theorem}[equation]{Theorem}
\newtheorem{lemma}[equation]{Lemma}
\newtheorem{proposition}[equation]{Proposition}
\newtheorem{corollary}[equation]{Corollary}

\theoremstyle{definition}
\newtheorem{example}[equation]{Example}
\newtheorem{remark}[equation]{Remark}

\newtheorem{bekezdes}[equation]{}
\newtheorem{definition}[equation]{Definition}

\newcommand{\gras}[1]{{\mathbb #1}}
\newcommand{\C}{\gras{C}}
\newcommand{\Z}{\gras{Z}}
\newcommand{\Q}{\gras{Q}}
\newcommand{\N}{\gras{N}}

\newcommand{\E}{\mathcal{E}}

\newcommand{\R}{\gras{R}}
\newcommand{\PP}{\gras{P}}

\newcommand{\xpq}{{\mathcal X}_{p,q}}
\newcommand{\calT}{{\mathcal T}}

\begin{document}
\title[Milnor fibers of cyclic quotient singularities]{On
the Milnor fibers of cyclic quotient
  singularities}
\author{\sc Andr{\'a}s N{\'e}methi}
\address{Renyi Institute of Mathematics, POB 127, H-1364 Budapest, Hungary.}
\email{nemethi@renyi.hu}
\author{\sc Patrick Popescu-Pampu}
\address{Univ. Paris 7 Denis Diderot, Inst. de
  Maths.-UMR CNRS 7586, {\'e}quipe "G{\'e}om{\'e}trie et dynamique" \\
  Site Chevaleret, Case
  7012, 75205 Paris Cedex 13, France.}
\email{ppopescu@math.jussieu.fr}



\date{April 30, 2009.}

%%% ----------------------------------------------------------------------
\begin{abstract}
  The oriented link of the cyclic quotient singularity
  $\mathcal{X}_{p,q}$ is orientation-preserving diffeomorphic to the
  lens space $L(p,q)$ and carries the standard contact structure
   $\xi_{st}$. Lisca classified  the Stein fillings of
  $(L(p,q), \xi_{st})$ up to diffeomorphisms
  and conjectured that they correspond
  bijectively through an {\it explicit} map to the Milnor fibers
  associated with the irreducible components (all of them being
  smoothing components) of the reduced
  miniversal space of deformations of $\mathcal{X}_{p,q}$. We prove
  this conjecture using the smoothing equations given by
  Christophersen and Stevens. Moreover, based on a different description
  of the Milnor fibers given by de Jong and van Straten,
  we also canonically identify these fibers with Lisca's fillings.
  Using these and a newly introduced additional structure --- the
  order --- associated  with lens spaces,
   we prove that the above Milnor fibers
  are pairwise non-diffeomorphic (by diffeomorphisms which preserve
  the orientation and order).
  This also implies that de Jong and van Straten parametrize in the same way
  the components of the reduced miniversal space of deformations
  as Christophersen and Stevens.
\end{abstract}

\maketitle
\pagestyle{myheadings} \markboth{{\normalsize
A. N{\'e}methi and P. Popescu-Pampu}}
{{\normalsize Milnor fibers of cyclic quotient singularities}}

% \tableofcontents

\section{Introduction} \label{intro}

\subsection{Lisca's conjecture}
In \cite{L 04}, Lisca announced a classification of the symplectic
fillings of the standard contact structure on
lens spaces up to orientation-preserving diffeomorphisms. Detailed
proofs were given in \cite{L 08}.

We recall briefly his classification. Let $L(p,q)$ be an oriented lens
space. Lisca provides first by surgery
diagrams a list of compact oriented $4$-manifolds
$W_{p,q}(\underline{k})$ with boundary $L(p,q)$. They
are parametrized by a set $K_r(\frac{p}{p-q})$ of
sequences of integers $\underline{k}\in \N^r$
(for its definition see (\ref{reprzerobis})).
He showed that each manifold $W_{p,q}(\underline{k})$
admits a  structure of Stein surface, filling the standard contact
structure on $L(p,q)$ and that any symplectic filling of this
standard contact structure is orientation-preserving diffeomorphic to
a manifold obtained from one of the $W_{p,q}(\underline{k})$ by a
composition of blow-ups (that is, in the language of differential
topology, by doing connected summing with the complex projective plane
endowed with the opposite orientation).

Particular cases of his theorem had been proved before by Eliashberg
\cite{E 90} (for $\gras{S}^3$) and McDuff \cite{M 90} (for the spaces
$L(p,1)$, for all $p \geq 2$).

In general, the oriented diffeomorphism type of the boundary and
the parameter $\underline{k}$ do not determine uniquely the
(orientation-preserving) diffeomorphism type of the fillings: for
some pairs the corresponding  types might coincide
(they are also listed by Lisca).

 Lisca noted that, following the works of
Christophersen  \cite{C 91} and Stevens \cite{S 91},
$K_r(\frac{p}{p-q})$ parametrizes also the irreducible components
of the reduced miniversal base space of deformations of the cyclic
quotient singularity $\mathcal{X}_{p,q}$. The oriented link of
this singularity is precisely a lens space $L(p,q)$. Each
component of the miniversal space is in this case a smoothing
component, that is, the generic local fiber over it is smooth. Its
oriented differentiable type is independent of the choice of the
generic point, and is called \emph{the Milnor fiber} of that
component. By construction, the Milnor fiber is orientation-preserving
diffeomorphic to a Stein filling of  ($L(p,q), \xi_{st})$. Lisca conjectured:

\bekezdes {\bf Conjecture} \cite[page 768]{L 08} {\em The Milnor fiber of
  the irreducible component of the reduced
miniversal base space of the cyclic quotient singularity
$\mathcal{X}_{p,q}$, parametrized in \cite{S 91} by
$\underline{k}\in K_{r}(\frac{p}{p-q})$,  is diffeomorphic to
$W_{p,q}(\underline{k})$.}

\vspace{3mm}

On the other hand, in \cite{JS 98}, de Jong and van Straten
studied by an approach completely different from Christophersen
and Stevens the deformation theory of cyclic quotient
singularities (as a particular case of sandwiched singularities).
They also parametrized  the Milnor fibers of $\mathcal{X}_{p,q}$
using the elements of the set $K_r(\frac{p}{p-q})$. Therefore, one
can formulate the previous conjecture for their parametrization as
well.

\subsection{}\label{ss:mainRe}  {\bf The main results} of
the present article and their consequences  are the following:

\medskip

$\bullet$  \emph{We introduce an additional structure associated
with any (non-necessarily oriented) lens space: the `order'. Its
meaning in short is the following: geometrically it is a (total)
order of the two solid tori separated by the (unique) splitting
torus of the lens space; in plumbing language, it is an order of
the two ends of the plumbing graph (provided that this graph has
at least two vertices). Then we show that the oriented
diffeomorphism type and the order of the boundary, together with
the parameter $\underline{k}$ determines uniquely the filling.}

\vspace{2mm}

%\noindent
$\bullet$  \emph{We endow in a natural way all the boundaries of
the spaces involved (Lisca's fillings $W_{p,q}(\underline{k})$,
Christophersen-Stevens' Milnor fibers $F_{p,q}(\underline{k})$,
and de Jong-van Straten's Milnor fibers
$W(\underline{a},\underline{k})$) with orders  --- the
corresponding spaces with these extra-structure will be
distinguished by  $^*$. Then we prove that all these spaces are
connected by orientation-preserving diffeomorphisms which preserve
the order of their boundaries: $W_{p,q}(\underline{k})^*\simeq
F_{p,q}(\underline{k})^*\simeq W(\underline{a},\underline{k})^*$.
This is an even stronger statement than the result expected by
Lisca's conjecture since it eliminates the ambiguities present in
Lisca's classification.}

\vspace{2mm}

% \noindent
$\bullet$  \emph{In fact, we even provide a fourth description of
the Milnor fibers: they are  constructed by a minimal sequence of blow ups
of the projective plane which eliminates the indeterminacies of a
rational function which depends on $\underline{k}$, see
(\ref{cor:9.2}) and (\ref{finalrem})(1). This is in the spirit of
Balke's work \cite{B 99}. }

\vspace{2mm}

  $\bullet$ \emph{As a byproduct it follows
(see \S \ref{ljs}) that both Christophersen-Stevens and de
  Jong-van Straten parame\-trized the components of the mini\-versal
  base space   in the same way (a fact not proved before, as far as we
  know).}

\vspace{2mm}

$\bullet$  \emph{Moreover, we obtain that the Milnor fibers
corresponding to the various irreducible components of the
miniversal space of deformations of $\mathcal{X}_{p,q}$ are
pairwise non-diffeomorphic by orientation-preserving
diffeomorphisms whose restrictions to the boundaries preserve the
order.}

\subsection{Symplectic fillings and singularities} Our work may be
considered as a continuation of the efforts to find all possible
Stein or, more generally, symplectic fillings of the contact links
of normal surface singularities. As a continuation of \cite{CP 04}, in
\cite{CNP 06} we showed with
Caubel that such contact structures are determined up to
contactomorphism by the topology of the link, that is, they depend
only on the topological type of the singularity, and not on its
analytical type. Therefore, singularity theory gives the following
Stein fillings up to diffeomorphism: the minimal resolution of
good representatives (which may be made Stein by deformation of
the complex structure, see \cite{BO 97}) and the Milnor fibers of
the smoothings of \emph{all} the analytical realizations of the given
topological type.

A natural question is \emph{to determine the topological types of
normal surface singularities for which one gets in this way
\emph{all} the Stein fillings of the associated contact manifold, up
to diffeomorphism}.

Ohta and Ono proved that this is the case for simple elliptic
singularities (see \cite{OO 03}) and for simple singularities,
that is, rational double points (see \cite{OO 05}). The above
positive answer to Lisca's conjecture shows that this is also the
case for cyclic quotient singularities.

We would like to stress some points regarding the previous classes
of singularities. Both simple and cyclic quotient singularities
are \emph{taut} singularities, that is, their analytical type is
determined by their topological type. Moreover, they are also
\emph{rational} singularities, 
hence their minimal resolution is diffeomorphic to the Milnor
fiber of one of the smoothing components, the so-called
\emph{Artin
  component} (see \cite[pages 33-34]{BR 95}). By contrast, simple
elliptic singularities are neither rational, nor taut, and their
minimal resolution is not diffeomorphic to the Milnor fiber of
some smoothing.

\subsection{Organisation of the paper}
In \S \ref{gener} we recall necessary facts about Hirzebruch-Jung
continued fractions and their geometric interpretation, while \S
\ref{cqls} contains some basic properties of cyclic quotient
singularities and lens spaces, expressed in terms of the geometry
of continued fractions. This section introduces the `order' of the
lens spaces too. Lisca's classification of the Stein fillings of
lens spaces is presented in  \S \ref{lisca}. In its last
subsection we reformulate his result using the notion of order. \S
\ref{chrste} contains the results of Christophersen and Stevens
regarding the structure of the reduced base of the miniversal
deformation of cyclic quotient singularities.
 In \S \ref{topsmooth} we recall de Jong and van
Straten's theory of deformations of sandwiched surface singularities
using decorated curves, which is specialized to cyclic quotient
singularities in \S \ref{cycl}.

In these preliminary six sections we provide several details on
the objects manipulated in order to try to make the paper readable
both by singularity theorists and topologists interested in
contact/symplectic topology. Considerable part of the preliminary
material is used in the proofs (and the remaining part is
conducive in the proper  understanding of the main
ideas/statements). On the other hand, even in these preliminary
sections, most of the `known' results are harmonized with the
newly introduced notion of order.

The main new results are contained in the last three sections.
 In \S \ref{lcs} we prove the `strong' version (cf.
subsection (\ref{ss:mainRe})) of Lisca's
conjecture using the equations of Christophersen and Stevens
describing the deformations of cyclic quotient singularities.

The identification of the Milnor fibers provided by the
construction of de Jong and van Straten with the Stein fillings is
done in \S \ref{ljs}. The proof needs a generalization of Lisca's
criterion for recognition of  each filling to a more homological
criterion, which is in turn proved in  \S \ref{capping}. The two
most important consequences are listed in \S \ref{concl}.

\subsection{Conventions and notations} \label{conv}
All the manifolds we consider will be oriented: any letter, say $W$,
 denoting a manifold will denote in fact an oriented  manifold.
We denote by $\overline{W}$ the manifold
obtained by changing the orientation of $W$,
and by $\partial W$ its boundary,
canonically oriented by the rule that the outward
normal \emph{followed} by the orientation of $\partial W$ gives the
orientation of $W$.

We work exclusively with (co)homology groups with
integral coefficients.

If $W$ is a  4-manifold, we denote by $Q_W:H_2(W)\times
H_2(W)\rightarrow \Z$ its intersection form
and by $\partial_W: H_2(W, \partial W)
\rightarrow H_1(W)$ the boundary
homomorphism. Additionally, if $W$ has a non-empty boundary, and
$S_1$ and $S_2$ are two 2-dimensional compact chains in $W$ with
disjoint boundaries which are  contained in $\partial W$, then
their intersection number is also well-defined and is denoted by
$S_1\cdot S_2$ or $S_1 \cdot_W S_2$.

If $a,b \in \N^*$ and $A$ is a commutative ring, we denote by
$\mbox{Mat}_{a,b}(A)$
the set of matrices with $a$ rows and $b$ columns with
coefficients in $A$.

If $M$ is an abelian group and $\gras{K}$ is a field, we write
$M_{\gras{K}}:=M \otimes_{\Z}\gras{K}$.

\subsection*{Acknowledgements}
The second author is grateful to Paolo Lisca for having attracted
his attention  on his conjecture. We are grateful to Jan Stevens
for acquainting us with the manuscript  \cite{B 99}
and to Duco van Straten for several helpful discussions.

\medskip
\section{Generalities on continued fractions
and duality of supplementary cones}
\label{gener}

\subsection{Hirzebruch-Jung continued fractions}\label{ss:3.1} If
$\underline{x}=(x_1,\ldots,x_n)$ are variables, the {\it
Hirzebruch-Jung continued fraction}:
\begin{equation} \label{contfrac}
  [x_1,\ldots,x_n]:= x_1 - \cfrac{1}{x_2 -
                          \cfrac{1}{\cdots - \cfrac{1}{x_n}}}
\end{equation}
can also be defined by induction on $n$ through the formulae:
$[x_1]=x_1$ and $[x_1,\ldots, x_n]=x_1-1/[x_2,\ldots,x_n]$ for $n
\geq 2$. One shows that:
\begin{equation} \label{contquot}
   [x_1,\ldots,x_n]= \frac{Z_n(x_1,\ldots,x_n)}{Z_{n-1}(x_2,\ldots,x_n)}
\end{equation}
where the polynomials $Z_n \in \Z[x_1,\ldots,x_n]$ satisfy
the inductive formulae:
\begin{equation} \label{recpol}
  Z_n(x_1,\ldots,x_n)=x_n\cdot
  Z_{n-1}(x_1,\ldots,x_{n-1})-Z_{n-2}(x_1,\ldots,x_{n-2}) 
  \mbox{ for all } n \geq 1,
\end{equation}
 with $Z_{-1} \equiv 0, \: Z_0 \equiv 1$ and $Z_1(x)=x$.
In fact, $Z_n(\underline{x})$ equals the determinant of the matrix
$M(\underline{x})\in \mbox{Mat}_{n,n}(\Z)$, whose entries are
$M_{i,i}=x_i$, $M_{i,j}=-1$ if $|i-j|=1$, $M_{i,j}=0$  otherwise.
Hence,  besides (\ref{recpol}), they satisfy many `determinantal
relations'. E.g.:
\begin{equation} \label{sym}
  Z_n(x_1,\ldots,x_n)=Z_n(x_n,\ldots,x_1).
\end{equation}

\begin{definition}  \cite{OW 77} \label{admis}
We say that  $\underline{x}\in \N^n$ is \emph{admissible} if
the matrix $M(\underline {x})$ is positive semi-definite
of rank $\geq n-1$.
 Denote by  $\mbox{adm}(\N^n)$ the set of admissible $n$-tuples.
\end{definition}

If $\underline{x} $ is admissible and  $n>1$, then each $x_i>0$.
Moreover, if $[x_1,\ldots,x_n]$ is admissible then
$[x_n,\ldots,x_1]$ is admissible too.

 Each rational
number $\lambda
>1$ admits a unique \emph{Hirzebruch-Jung continued fraction
expansion} (in short, a \emph{HJ-expansion}) of the form:
$$\lambda=[x_1,\ldots,x_n], \ \mbox{where $\: x_i
                        \in \N, \:
x_i \geq 2$ \ for all \ $i \in \{1,\ldots,n\}$}.$$

\subsection{The geometrical interpretation.}\label{3:2}
 Next we explain an interpretation of the HJ-expansions using affine
 geometry, see e.g. \cite{Oda} and \cite{PPP 07}.

Consider a free abelian group $N$ of rank two.
   An {\it oriented cone} is a rational strictly convex cone
   $\sigma$ in $N_{\R}$ with a choice of an order of its edges.
   We denote the two primitive elements of $N$ which
   generate the edges by $e_1$ and $e_2$, where $(e_1,e_2)$ is the
   order of the edges.
    We denote by $\sigma'$ the same cone with the opposite choice
    $(e_2,e_1)$ of order
   of its edges and by $\overline{\sigma}$ the {\it supplementary
     cone} generated by $(-e_1, e_2)$.

Consider the convex hull of the points of $N$ situated in
$\sigma\setminus 0$ and denote by $P_{\sigma}$ the union of the
compact edges of its boundary. It is a finite polygonal line
joining $e_1$ to $e_2$. Denote by $v_0, \ldots, v_{s+1}$ the
lattice points situated on it, in the order in which they are
encountered when one travels from $v_0:= e_1$ to $v_{s+1}:=e_2$.
Then for some integers $b_i\geq 2$:
\begin{equation} \label{basrel}
    v_{i-1} + v_{i+1}= b_i v_i
\ \ \mbox{ for all } \ i \in \{1,\ldots,s\}.
\end{equation}

Write also $e_2=-q v_0 +p v_1$. Then  $p$ and $q$ are coprime with
$p>q >0$ provided that $(e_1,e_2)$ is not a basis of $N$, and
$p=1$, $q=0$ otherwise. From now one we suppose that we are in the
first case.

\begin{lemma}\label{lem:2.2.3}
   With the previous notations,  $\frac{p}{q}=[b_1,\ldots,b_s]$.
\end{lemma}

Both $\frac{p}{q}$ and the sequence $(b_1,\ldots,b_s)$ are
complete invariants of the pair $(N, \sigma)$, up to isomorphisms
(that is, isomorphisms of free groups which send one cone onto the
other and preserve the order of the edges).  We say that
$\frac{p}{q}$ is \emph{the type} and $(b_1,\ldots,b_s)$ \emph{the
associated sequence} of the oriented cone $(N, \sigma)$. If one
changes the orientation of the cone (i.e. $\sigma$ into
$\sigma'$), then the type of $(N, \sigma')$ becomes
$\frac{p}{q'}$, where $q'$ is the unique positive number such that
$q'<p$ and $qq'\equiv 1 \: (\mbox{mod} \: p\:)$; the associated
sequence becomes $(b_s,\ldots,b_1)$.

\begin{figure}[ht]
\labellist \small\hair 2pt \pinlabel {$0$} at 253 -10 \pinlabel
{$e_1$} at 72 -10 \pinlabel {$-e_1$} at 433 -10 \pinlabel {$v_0$}
at 52 20 \pinlabel {$\overline{v}_0$} at 451 20 \pinlabel {$v_1 $}
at 99 27 \pinlabel {$\overline{v}_1$} at 410 27 \pinlabel {$v_2$}
at 152 44 \pinlabel {$\overline{v}_2$} at 373 47 \pinlabel {$v_s$}
at 285 210 \pinlabel {$\overline{v}_r$} at 371 214 \pinlabel
{$v_{s+1}=e_2=\overline{v}_{r+1}$} at 336 276 \pinlabel
{$\sigma\simeq \sigma_{p,q}$} at 84 189 \pinlabel {$\mbox{here one
reads}$} at 84 159 \pinlabel {$(b_1,\ldots,b_s)$} at 84 139
\pinlabel {$\overline{\sigma}\simeq \sigma_{p,p-q}$} at 492 189
\pinlabel {$\mbox{here one reads}$} at  492 159 \pinlabel
{$(a_1,\ldots,a_r)$} at 492 139 \pinlabel {$P_{\sigma}$} at 105 86
\pinlabel {$P_{\overline{\sigma}}$} at 472 86
\endlabellist
\centering
\includegraphics[scale=0.50]{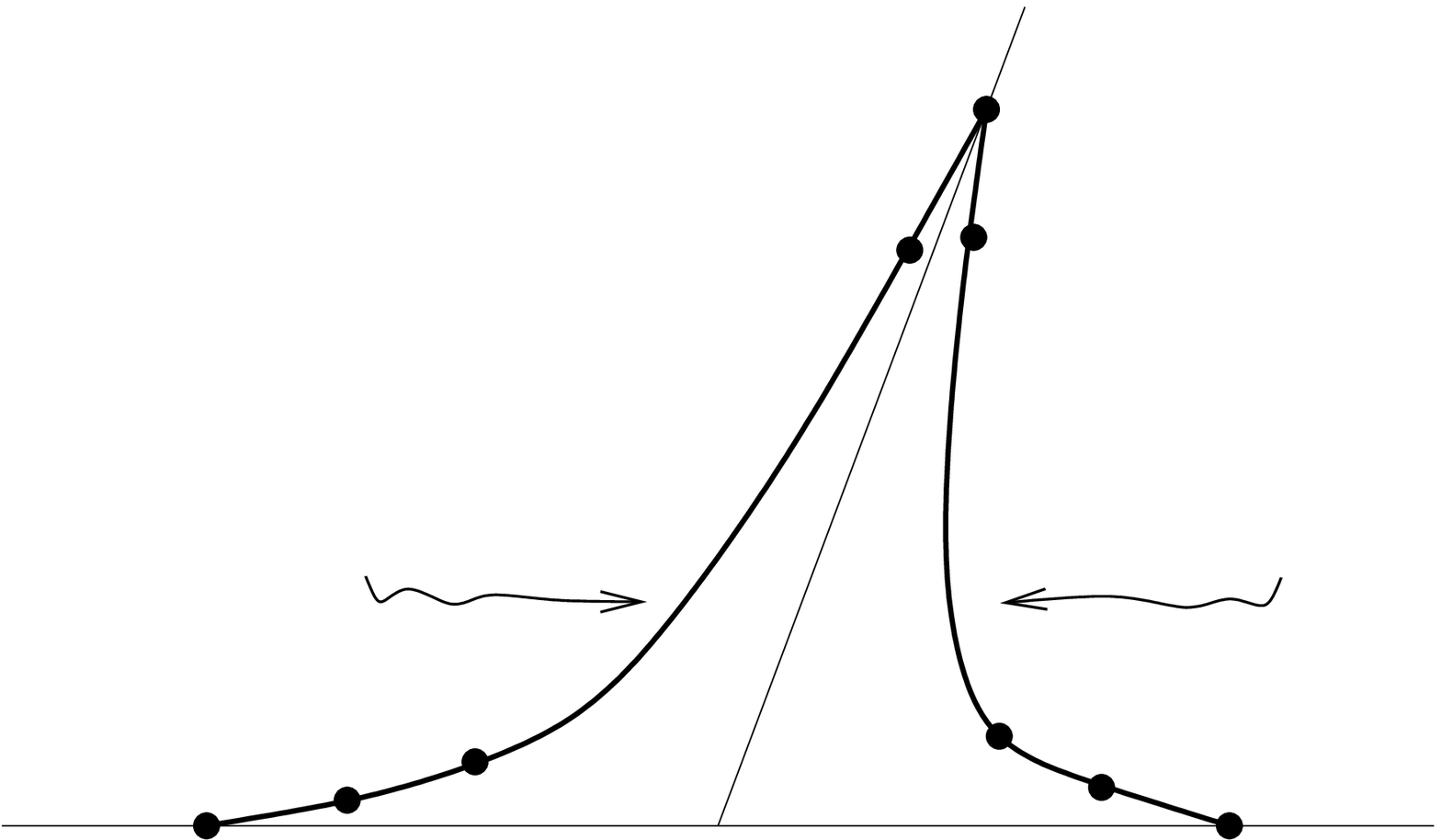}
\caption{Two supplementary cones}
\label{fig:Suppl}
\end{figure}

Consider now both
HJ-expansions:
\begin{equation} \label{dualexp}
  % \left\{ \begin{array}{l}
       \displaystyle{\frac{p}{q}}= [b_1,\ldots,b_s]
\ \ \  \mbox{and} \ \ \
       \displaystyle{\frac{p}{p-q}}= [a_1,\ldots,a_r].
  %           \end{array}  \right.
\end{equation}
There is a \emph{duality} of the sequences $\underline{a}$ and
$\underline{b}$. Its geometric interpretation is the following.
Start from an oriented cone $\sigma\simeq \sigma_{p,q}$ in
$N_{\R}$ of type $\frac{p}{q}$ (well-defined up to \emph{unique}
isomorphism). Consider its supplementary cone $\overline{\sigma}$,
the polygonal line $P_{\overline{\sigma}}$ and the sequence
$(\overline{v}_0,\ldots, \overline{v}_{r+1})$ of lattice points on
it, starting from $\overline{v}_0:= -e_1$ and ending with
$\overline{v}_{r+1}:= e_2$ (see Figure \ref{fig:Suppl}).

\begin{lemma} \label{dualseq}
   $\overline{\sigma}\simeq
  \sigma_{p,p-q}$, i.e. the sequence associated with
  $(N, \overline{\sigma})$ is $(a_1,\ldots,a_r)$.
\end{lemma}

There is an important point we wish to emphasize. In general, in
torical presentations,  cones and their dual cones are sitting in
different  lattices (dual to each other).  Here, though
 the supplementary cone can canonically be identified with
the dual cone (using the area-symplectic form equal to $1$ on a
positive basis of the lattice,  see \cite[page 145]{PPP 07}),  both
$\sigma$ and 
$\overline{\sigma}$  are represented in the {\it same} lattice
$N$. This will allow us to connect by linear relations vectors
from both cones in the same $N$ (see Theorem (\ref{precdual})). Such a
relation interpreted in 
homology is important in the proof of the main result from \S
\ref{ljs}.

There is a canonical way to rewrite the entries of the continued
fractions, where $(2)^\ell$ means the constant sequence with
$\ell$ terms all equal to $2$, and $n_j \geq 3 $ for all $1\leq
j\leq t$:
\begin{equation} \label{HJexp1}
    \frac{p}{p-q}=[(2)^{m_1-1}, n_1, (2)^{m_2-1},n_2,\ldots, n_t,
    (2)^{m_{t+1}-1}].
\end{equation}
Geometrically, $t\in\N$ is the number of `interior' vertices of the
polygonal line $P_{\overline{\sigma}}$, $(m_1,...,m_{t+1})$ the
sequence of integral lengths of the edges of
$P_{\overline{\sigma}}$ (hence $ m_i \geq 1$ for all $1\leq i
\leq  t+1$).

From the arithmetical point of view of the HJ-expansions of $p/q$ and
$p/(p-q)$, the duality is reflected by
Riemenschneider's point diagram \cite{R 74} which basically  says
that:
\begin{equation} \label{HJexp2}
   \frac{p}{q}=[m_1+1, (2)^{n_1-3}, m_2 +2,(2)^{n_2-3},\ldots, m_t+2,
   (2)^{n_t-3}, m_{t+1}+1].
\end{equation}
In particular:
\begin{equation}\label{eq:S4}
r=1+\sum_{1\leq i\leq s}(b_i-2)=-1+\sum_{1\leq i\leq t+1}m_i \ \ \
\mbox{and} \ \ \ s=1+\sum_{1\leq j\leq t}(n_j-2).
\end{equation}
But there is an even deeper relation at the level of $N$ (see
\cite[Proposition 5.3]{PPP 07}):
\begin{theorem} \label{precdual}
   Set
  $w_l:= v_{1 + \sum_{1 \leq j \leq l-1} (n_j -2)}$  \, for all \,
   $1\leq l\leq t+1$.
  Then:
  \begin{equation} \label{fundiff}
    \overline{v}_{i+1}-\overline{v}_i =w_l \  \mbox{ if } \ m_1 + \cdots
    m_{l-1} \leq i \leq m_1 + \cdots + m_l -1.
  \end{equation}
\end{theorem}
For a detailed discussion of similar relations connecting a cone
with its supplementary cone, and their relationship with continued
fractions, see \cite{PPP 07}.

\medskip
\section{Generalities on
cyclic quotient singularities
and lens spaces}
\label{cqls}

We recall the definitions of \emph{cyclic quotient
  singularities} and  \emph{lens spaces}. Additionally,  we
 introduce the notion of {\em order} associated with a lens space
 and we discuss its relationship
 with the group of automorphisms and dualities.
 See \cite[pages 99--105]{BHPV 04}  for a  classical presentation
of cyclic quotient singularities, \cite{N 81,PPP 07} for details about
plumbings of links of  surface singularities, and   \cite{B 83,PPP 07}
for  the geometry of the splitting torus of lens spaces.

\subsection{The definitions}\label{ss:4.1}
 Let $p, q$ be coprime  integers such that $p>q>0$.

\begin{definition} \label{defhj}
  The {\it cyclic quotient} (or   {\it Hirzebruch-Jung})  {\it singularity}
  $({\mathcal X}_{p,q},0)$   is the germ
  of the quotient  ${\mathcal X}_{p,q}$ of $\C^2$ by the action
  $\xi*(x,y)=(\xi x, \xi^q y)$
  of the cyclic group $\{\xi \in \C\:,\: \xi^p=1\}\simeq \Z/p\Z$.
  Its oriented link is the (oriented)    {\it lens space} $L(p,q)$.
\end{definition}

In particular,  $L(p,q)$ is  the quotient of $\gras{S}^3$ by the
above action of $\Z/p \Z$  (this definition does not include
$\gras{S}^3$ and $\gras{S}^2 \times \gras{S}^1$, which are
sometimes also considered to be lens spaces).  Bonahon in \cite{B
83} proved that {\em each lens space contains an embedded
2-dimensional torus
--- a so-called splitting torus ---, unique up to an isotopy,
which bounds on each side a solid torus}. The set $\calT$ of solid
tori bounded by a splitting torus, identified modulo isotopies of
the ambient space, is a set of one or two elements. It has  one
element exactly when the solid tori can be interchanged by an
isotopy. This happens exactly when $q\in\{1,p-1\}$, cf. \cite[page
308]{B 83}.

Let us define an additional structure associated with an
(unoriented) lens space. It has a similar nature as the notion of
orientation ($\calT$ is analogous to the set of connected
components of the orientation bundle of a manifold), but it is
independent of it.

\begin{definition}\label{def:order} An {\em order} of an (unoriented) lens
  space  is a total  order on the set $\calT$.
\end{definition}
Clearly, if $q\in\{1,p-1\}$, then this supports no additional
information.  In all other cases the order
 distinguishes \emph{the first} and \emph{the second}  of the two
 (non-isotopic) solid tori bounded by \emph{any} splitting torus.

The unicity of the splitting torus $\tau$ allows one to associate
with any (unoriented) lens space $L$ a free abelian group of rank
2, namely  $N:= H_1(\tau)$. In the sequel $\nu_*:N\to H_1(L)$ will
stay for the homological morphism induced by the inclusion
$\nu:\tau\hookrightarrow L$.

\begin{remark}\label{remark:3.1}
In fact, $N$ is well-defined up to the induced action of the
isotopies of the lens space which move $\tau$ into itself. This is
non-trivial only if $q\in\{1,p-1\}$, but even in those cases any
such induced isomorphism $\varphi_N: N\to N$ satisfies $\nu_*\circ
\varphi_N=\nu_*$.
\end{remark}

\subsection{The order and its type}\label{ss:4.2}
We explain now a way to extract the numbers
 $\{\frac{p}{q},\frac{p}{q'}\}$
from an oriented lens space diffeomorphic to $L(p,q)$.

Let $L$ be an \emph{oriented} lens space. Choose an order of the
two solid tori bounded by a splitting torus $\tau$: denote by
$L_1$ the first  and by $L_2$ the second one. Orient $\tau$ as the
boundary of $L_1$. Therefore, $N=H_1(\tau)$ gets an induced
orientation (dual to the orientation of $H^1(\tau)$ such that the
cup product of a positive basis is positive on the fundamental
class of $\tau$). There is up to isotopy a unique \emph{meridian}
of $L_1$ (that is, an unoriented simple closed curve on $\tau$
which is non-trivial homotopically on $\tau$ but is trivial in
$L_1$). Orient it in an arbitrary way and denote by $e_1\in N$ its
homology class. Then orient the meridian of $L_2$ such that its
homology class $e_2 \in N$ forms a positive basis $(e_1, e_2)$ of
$N$ with respect to the orientation defined before. Denote by
$\sigma$ the oriented strictly convex cone in $N_{\R}$ generated
by $e_1$ and $e_2$, taken in this order. Let $\frac{p}{q}$ be the
type of the oriented cone $(N, \sigma)$, cf. (\ref{3:2}).
\begin{lemma}
  $L$ is (orientation-preserving) diffeomorphic to $L(p,q)$.
\end{lemma}

By choosing the opposite orientation of the meridian of $L_1$, one
gets $(-e_1,-e_2)$ instead of $(e_1,e_2)$, hence $-\sigma$ instead
of $\sigma$, whose type is the same $\frac{p}{q}$.
By changing the order of the solid tori, one gets
as new type $\frac{p}{q'}$. This  also reproves the classical fact
that $L(p,q)$ is orientation-preserving diffeomorphic to $L(p, q_1)$
if and only if $q_1\in\{q,q'\}$.

If $\#\calT=1$, then $\frac{p}{q}=\frac{p}{q'}$. If $\#\calT=2$
and we fix an order, then  we get without ambiguity a unique
element of  $\{\frac{p}{q},\frac{p}{q'}\}$. Hence, an order always
provides a well-defined  element of the set
$\{\frac{p}{q},\frac{p}{q'}\}$,  called the {\em type of the
order}.

If $q'\not=q$ then from the type of the order one can recover the
order itself. Indeed, the type contains the information regarding
the oriented cone, whose ordered edges correspond to an order of
the two meridians. This is not the case for
$q'=q\not\in\{1,p-1\}$, since $\#\calT=2$ but
$\#\{\frac{p}{q},\frac{p}{q'}\}=1$ (hence the order is a sharper
invariant than its type).

 Notice that
if we change the orientation of the above lens space $L(p,q)$,
then in the above construction $(e_1,e_2)$ can be replaced by
$(-e_1,e_2)$ (i.e. $\sigma$  by $\overline{\sigma}$ of type
$\frac{p}{p-q}$), hence $\overline{L(p,q)}=L(p,p-q)$.

\begin{remark}
The notation `$L(p,q)$' is not uniform in the literature:
sometimes, what we call $L(p,q)$ is denoted by $L(p,p-q)$; the
ambiguity originates in the orientation choice.
\end{remark}

\subsection{Self-diffeomorphisms}\label{ss:4.4}
The previous discussion  also provides the group of
orien\-tation-preserving self-diffeomorphisms Diff$^+(L)$ of an
oriented lens space $L$, cf.  \cite{B 83}. For any $\varphi\in
\mbox{Diff}^+(L)$ write $\varphi_*$ for the induced morphism at
the level of $H_1(L)\simeq \Z/p\Z$. Then the isotopy class of
$\varphi$ is uniquely determined by $\varphi_*$, and $\varphi_*$
can only be multiplication by $\pm 1$ or $\pm q$, and $\pm q$
occurs only if $q'=q$ (corresponding to how $\varphi$  changes the
orientation of the meridians and/or the  solid tori).

If $q'\not= q$, then Diff$^+(L)$ preserves automatically the order
(that is, any of the two possible orders is left invariant). If
$q=q'\notin \{ 1, p-1\}$, then $\varphi$ reverses it exactly when
$\varphi_*$ is the multiplication by $\pm q $. 
Hence, once an order $o$ is fixed, the subgroup Diff$^{+,o}(L)$ of
Diff$^+(L)$ which preserves $o$ is isomorphic to $\Z/2\Z$, and its
elements induce $\varphi_*=\pm1$ on $H_1(L)\simeq \Z/p\Z$, with the
exception of $p = 2$, when Diff$^{+,o}(L)\simeq \{\mathrm{Id}_{\Z/2\Z}\}$.

\subsection{The order and the sequence
$\{v_1,\ldots,v_s\}$}\label{ss:4.new} If both an  orientation and
an order are fixed on a lens space (say by the choice of the
numbering $L_1$ and $L_2$ of the solid tori bounded by a splitting
torus), then the discussion from (\ref{ss:4.2}) shows that the
sequence $(v_1,\ldots,v_s)$ is well-defined up to a sign and the
automorphisms $\varphi_N$ of $N$ from (\ref{remark:3.1}).

For such a situation later we will use the following (ordered) set
of elements of $H_1(L)$, well-defined up to a sign (cf.
 (\ref{remark:3.1})), associated with an oriented and ordered $L$:
\begin{equation} \label{alphas}
   (\alpha^L_1,\ldots, \alpha^L_s):= \pm (\nu_*(v_1),\ldots,\nu_*(v_s))
\in H_1(L)^s.
\end{equation}
If  in the above construction one changes the orientation of $L$
and  one keeps {\em the same order}, then one gets:
\begin{equation} \label{alphabars}
   (\alpha^{\overline{L}}_1,\ldots, \alpha^{\overline{L}}_r):=
  \pm  (\nu_*(\overline{v}_1),\ldots,\nu_*(\overline{v}_r))\in H_1(L)^r,
\end{equation}
where the $v_i$ and $\overline{v}_i$ are related as in subsection
(\ref{3:2}).

\subsection{Order and plumbings}\label{ss:4.3}
Consider again the surface $\mathcal{X}_{p,q}$. The sequences
$\underline{a}$ and $\underline{b}$ from (\ref{dualexp}) guide two
different geometrical packages: $\underline{a}$ is related more to
the equations  and deformations of $\mathcal{X}_{p,q}$, see e.g.
\S \ref{chrste} (\cite{C 91} and \cite{JS 98} use even the
notation $\mathcal{X}_{p,q}=X(a_1,\ldots,a_r)$), while
$\underline{b}$ is related to the resolution.

\begin{figure}[ht!]

\begin{picture}(300,35)(-30,0)
\put(20,10){\circle*{4}} \put(60,10){\circle*{4}}
\put(160,10){\circle*{4}} \put(200,10){\circle*{4}}
\put(20,10){\line(1,0){60}} \put(140,10){\line(1,0){60}}
\put(20,25){\makebox(0,0){$-b_1$}}
\put(60,25){\makebox(0,0){$-b_2$}}
\put(160,25){\makebox(0,0){$-b_{s-1}$}}
\put(200,25){\makebox(0,0){$-b_s$}}
\put(110,10){\makebox(0,0){$\ldots$}}
\end{picture}

\caption{The  graph  $G(\underline{b})$} \label{fig:Gb}
\end{figure}
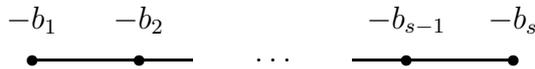

The dual graph $G(\underline{b})$  of the minimal resolution is  a
string, see Figure \ref{fig:Gb}. Let $\Pi(\underline{b})$ be the
oriented compact $4$-manifold with boundary obtained by plumbing
according to $G(\underline{b})$ (it is diffeomorphic to the
minimal resolution of a Milnor representative of
$(\mathcal{X}_{p,q},0)$). It contains oriented 2-spheres
$\{S_i\}_{1 \leq i  \leq s}$ which intersect according to the
graph $G(\underline{b})$ and which are realized in the minimal
resolution  by the irreducible exceptional curves. By construction
$\partial \Pi(\underline{b})=L(p,q)= \overline{\partial
\Pi(\underline{a})}$.

Notice also that the permutation $(x,y)\mapsto (y,x)$ of the
coordinates of $\C^2$ realizes an isomorphism ${\mathcal
X}_{p,q}\to {\mathcal X}_{p,q'}$. Hence, a priori there is no
preferred order of the coordinate axes or of the 2-spheres
$\{S_i\}_{1 \leq i  \leq s}$ if $s>1$:
 the marking  $\{S_i\}_{1\leq i  \leq s}$ can be replaced freely by
 $\{S_{s-i}\}_{1 \leq i  \leq s}$.
But the ambiguity disappears (for both plumbing graphs
$G(\underline{a})$ and $G(\underline{b})$) if $\#\calT=2$ and an
order is fixed. We explain this in the next paragraphs.

With any marking $\{S_i\}_{1\leq i\leq s}$ of the 2-spheres,
associate a collection $(R_i)_{1 \leq i
  \leq s}$ of pairwise disjoint
  oriented discs properly embedded inside $\Pi(\underline{b})$, $R_i$
   being a normal slice of $S_i$ inside the associated part
  of the plumbing decomposition of $\Pi(\underline{b})$, with $R_i
  \cdot_{\Pi(\underline{b})} S_i =+1$.

Assume that $q\not\in\{1,p-1\}$, or equivalently,
$s\not=1\not=r$, i.e. the lengths of both  $\underline{a}$  and
$\underline{b}$ are greater than one. Let us concentrate now on
$G(\underline{b})$ (there is a symmetric discussion for
$G(\underline{a})$ too). Any edge of the graph, via the plumbing
construction, determines a splitting torus. Then $\partial R_1$
and $\partial R_s$ sit in two different (non-isotopic) solid tori.

Each edge of $G(\underline{b})$ corresponds to a splitting torus
of $L$. Choose an ordering  $(L_1,L_2)$ of the solid tori bounded
by it such that the type of $L$ endowed with the associated order
is $p/q$. Then we mark the 2-spheres in such a way that
 $\partial R_1\subset L_1$ and $\partial R_s\subset L_2$.
This defines an order of 2-spheres. By this convention, not only
the type $p/q$ of the order and the ordered sequence $(b_1,\ldots
,b_s)$ of  weights of the graph correspond by the expansion
(\ref{lem:2.2.3}), but even if this sequence is symmetric, we
indicate  in the plumbed manifold $\Pi(\underline{b})$ which
2-sphere has index 1, respectively $s$. In fact, this
specification is equivalent with a choice of an order.

(Regarding Figure \ref{fig:Gb}, notice the following. If $s>1$,
and we replace the `symbols' $b_i$ by some integers $\geq 2$, then
we get an unordered graph. On the other hand, with the present
decoration we indicate which end-edge is $S_1$ respectively $S_s$,
hence Figure  \ref{fig:Gb} in fact represents an ordered graph
providing an order of its plumbed lens space.)

All the constructions of the present article are guided by
plumbing graphs. Using this, we will introduce uniformly an order
in all the lens spaces involved.

Assume  $s>1$ and $r>1$.  Fix the graph $G(\underline{b})$ and its
marking as in Figure  \ref{fig:Gb}. It defines an order (via the
rule $R_1\subset L_1$) of the plumbed manifold $\partial
\Pi(\underline{b})$.

Construct compatible markings on $G(\underline{a})$ and
$G(\underline{b})$ together with a canonical identification of the
plumbing manifolds $\partial \Pi(\underline{b})$ and $\partial
\Pi(\underline{a})$ as follows. Fix $G(\underline{b})$ with its
marking, which provides an order of $\Pi(\underline{b})$. The
orientation change of $\Pi(\underline{b})$ corresponds to the
change of each weight $-b_i$ into $b_i$ (but keeping the 1 resp.
$s$ marking of the ends). Notice that by a plumbing calculus of
the plumbing graphs we may keep track of the order of the lens
space (the order of the end-vertices  of the graph). Hence, when
we proceed a plumbing calculus in order to replace the graph
marked with $b_i$'s into one which has weights $\underline{a}$,
there is only one way to order/mark the sequence $\underline{a}$
in such a way that the order of the lens space will reflect
properly the order of the sequence $\underline{a}$.

Finally,  if a plumbing graph contains a subgraph having the
information regarding the above markings of $G(\underline{a})$ or
$G(\underline{b})$, we mark it in a compatible way.

\begin{definition}\label{preferred}
The order fixed in this way by the marked graphs is called the
{\it preferred order} of the lens spaces. $L$ endowed with a
preferred order will be denoted by $L^*$.
\end{definition}

\subsection{Order and coordinates.}\label{markcoord}
 An order of the plumbing/resolution graph can be
related with a choice of an order of the coordinates of $\C^2$ in
the definition (\ref{defhj}) as follows.

Assume first that $\underline{b}$ is not symmetric. If we fix $p$,
$q$ and the coordinates $(x,y)$ as in (\ref{defhj}), and we
mark/order the 2-spheres of the resolution in such a way that the
expansion $p/q=[b_1,\ldots,b_s]$ holds, then the strict transform
of (the image) of the $x$-axis (respectively, of $\{x=0\}$) can be
isotoped to $R_1$ (respectively $R_s$). We take this as a general
compatibility property even if $\underline{b}$ is symmetric: the
(order of the) coordinates $(x,y)$ is compatible with the ordering
of the 2-spheres if the above fact holds.

\subsection{} Finally, we verify another compatibility property. Notice that
each  $\partial R_i$ can be isotoped in  $\partial
\Pi(\underline{b})$  in a tubular neighborhood of $\tau$  (whose
first homology is  $N$).
\begin{proposition} \label{repres} \cite[page 176-177]{PPP 07}
Let $L$ be an ordered and oriented lens space.
If $q\not\in\{1,p-1\}$, then
  $([\partial R_1],\ldots, [\partial R_s])=
  \pm (v_1,\ldots,v_s)$   in $N$. In fact, in all cases,
  $(\nu_*[\partial R_1],\ldots, \nu_*[\partial R_s])=
  \pm (\nu_*(v_1),\ldots,\nu_*(v_s))$   in $H_1(L)$.
\end{proposition}

\medskip
\section{The Stein fillings of lens spaces,
following Lisca}\label{lisca}

The field of complex lines tangent to $\gras{S}^3$
is left invariant by the cyclic action used in
(\ref{defhj}), hence it descends to the so-called {\em standard
 contact structure} $\xi_{st}$ on $L(p,q)$.

\subsection{Lisca's construction}\label{ss:5.2} The next set is the
parameter set in the three main constructions presented in the body of
the paper.
\begin{definition}
For $r\geq 1$, denote by:
\begin{equation} \label{reprzero}
   K_r:=\{ \underline{k}= (k_1,\ldots ,k_r)\in  \mbox{adm}(\N^r) \ \
   |\   \ [k_1,\ldots ,k_r]   =0\}
\end{equation}
the set of \emph{admissible sequences which represent $0$}.
  For $\underline{k}=(k_1,\ldots ,k_r)\in K_r$  set
  $\underline{k}':=(k_r,\ldots,k_1)\in K_r$.
\end{definition}

We wish to emphasize  that the condition of \emph{admissibility} (cf.
\ref{admis}) is really restrictive. For example,
$\underline{k}=(2,1,1,1,1,2)\notin K_6$ although
$[\underline{k}]=0$. By admissibility, if $r>1$, then each
$k_i>0$. $K_1$ has only one element, namely $(0)$.

\begin{figure}[h!]
\labellist \small\hair 2pt 
\pinlabel{$\mathcal{L}(\underline{k})$} at -40 150
\pinlabel{$\mathbf{L(\underline{a},
\underline{k})}$} at 10 80 \pinlabel  {$a_1 -k_1$} at 130 -10
\pinlabel {$a_2-k_2$} at 237 -10 \pinlabel {$a_{r-1}-k_{r-1}$} at
505 -10 \pinlabel {$a_r-k_r$} at 621 -10 \pinlabel  {$k_1$} at 125
230 \pinlabel  {$k_2$} at 240 230 \pinlabel  {$k_{r-1}$} at 510
230 \pinlabel  {$k_r$} at 620 230 \pinlabel  {$-1$} at 95 35
\pinlabel  {$-1$} at 150 35 \pinlabel  {$-1$} at 210 35 \pinlabel
{$-1$} at 260 35 \pinlabel  {$-1$} at 485 35 \pinlabel  {$-1$} at
530 35 \pinlabel  {$-1$} at 590 35 \pinlabel  {$-1$} at 640 35
\pinlabel  {$\dots$} at 380 130 \pinlabel  {$\dots$} at 380 170
\pinlabel  {$\dots$} at 134 52 \pinlabel  {$\dots$} at 241 52
\pinlabel  {$\dots$} at 511 52 \pinlabel  {$\dots$} at 623 52
\endlabellist
\centering
\includegraphics[scale=0.50]{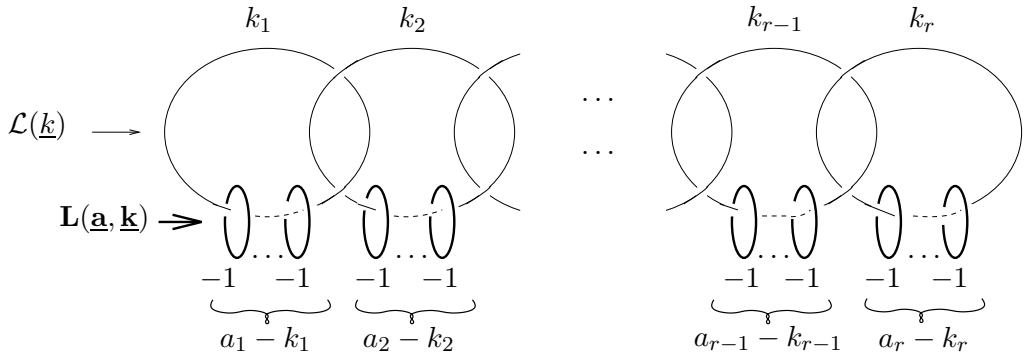}
\vspace{2mm} \caption{The framed link $L(\underline{a},
  \underline{k}) \subset N(\underline{k})$}
\label{fig:Surgery2}
\end{figure}

For two coprime integers $p,q$ with $p>q>0$ and HJ-expansion
$\frac{p}{p-q}=[a_1,\ldots,a_r]$, set:
\begin{equation} \label{reprzerobis}
   K_r(\textstyle{\frac{p}{p-q}})=K_r(\underline{a}):=\{ \underline{k}\in
   K_r \: | \: \underline{k}\leq \underline{a}\} \subset K_r.
\end{equation}
(Here, $\underline{k}\leq \underline{a}$ means that $k_i\leq a_i$
for all $i$.) Fix  an element $\underline{k} \in
K_r(\underline{a})$. Let ${\mathcal L}(\underline{k})$ be the
framed link of Figure \ref{fig:Surgery2} with $r$ components and
decorations $k_1,\ldots, k_r$ (i.e. the thick components are
neglected for a moment). Let $N(\underline{k})$ be the closed,
oriented 3-manifold given by surgery on $\gras{S}^3$ along
${\mathcal L}(\underline{k})$. Using the slam-dunk operation on
rationally-framed links in $\gras{S}^3$ (see \cite[page 163]{GS
99}), one obtains an orientation-preserving diffeomorphism:
\begin{equation} \label{diffsd}
\eta:  N(\underline{k}) \longrightarrow \gras{S}^1 \times
  \gras{S}^2.
\end{equation}

\begin{definition} \cite[page 766]{L 08}  \label{defStein}
  Consider the diffeomorphism $\eta$ from (\ref{diffsd}) and denote by
  $L(\underline{a}, \underline{k})\subset N(\underline{k})$  the
  thick framed link drawn in Figure \ref{fig:Surgery2}. Define
  $W_{p,q}(\underline{k})$ to be the smooth oriented 4-manifold with
  boundary obtained by attaching two-handles to $\gras{S}^1 \times
  \gras{D}^3$ along the framed link $\eta(L(\underline{a},
  \underline{k}))\subset \gras{S}^1 \times \gras{S}^2$.
\end{definition}
From the main theorem (1.1) of \cite{L 08}, one can extract the
following:
\begin{theorem}\label{th:lis} (a)
   All the manifolds $W_{p,q}(\underline{k})$ admit Stein structures
   which fill  \linebreak
    $(L(p,q), \xi_{st})$, and any Stein filling of $(L(p,q), \xi_{st})$ is
    diffeomorphic  to one of the manifolds
    $W_{p,q}(\underline{k})$.

    (b) $W_{p,q_1}(\underline{k}_1)$ is orientation-preserving
    diffeomorphic to
    $W_{p,q_2}(\underline{k}_2)$  if and only if
    $(q_2,\underline{k}_2)=(q_1,\underline{k}_1)$ or
     $(q_2,\underline{k}_2)=(q'_1,\underline{k}'_1)$.
\end{theorem}

\medskip

\subsection{Lisca's criterion to recognize
$W_{p,q}(\underline{k})$.}\label{ss:5.3} Once (\ref{th:lis})(a) is
proved, (\ref{th:lis})(b) follows from a straightforward
homological computation, whose essence is highlighted by the next
criterion. Let $W$ be a Stein filling of $(L(p,q), \xi_{st})$. Let
$V$ be the closed 4-manifold obtained by  gluing $W$ and
$\Pi(\underline{a})$ via an orientation-preserving diffeomorphism
$\phi: \partial W \rightarrow \partial
\overline{\Pi(\underline{a})}$ of their boundaries. Let $\mu:
\Pi(\underline{a}) \hookrightarrow V$ be the inclusion morphism.

 \begin{proposition} \label{glueplumb}\cite[\S 7]{L 08}
  Denote by $\{s_i\}_{1 \leq i \leq r}$ the classes of
     2-spheres $\{S_i\}_{1 \leq i \leq r}$ in
     $H_2(\Pi(\underline{a}))$ (listed in the same order as
     $\{a_i\}_{1\leq i\leq r}$), and also their images via the
     monomorphism $\mu_*:H_2(\Pi(\underline{a}))\rightarrow H_2(V)$.
      Then there is a  $\underline{k}\in K_r(\underline{a})$
      such that
    \begin{equation}\label{eq:ki}
\begin{small}
    \#\{  e \in H_2(V) \: | \ \ e^2=-1,\  s_i \cdot e\neq 0,\ s_j\cdot e
    =0 \ \mbox{for all} \
       j \neq i \}= 2 (a_i-k_i)
         \end{small}\end{equation}
         for all $i \in \{1,\ldots,r\}$. In this way one
    gets the pair $(\underline{a},\underline{k})$,  and
    $W$ is orientation-preserving diffeomorphic to
    $W_{p,q}(\underline{k})$.
 \end{proposition}
Notice that, as  $\{S_i\}_{1\leq i\leq r}$ and $\{S_{r-i}\}_{1\leq
i\leq r}$ cannot be  distinguished, the algorithm does not
differentiate $(\underline{a},\underline{k})$ from
$(\underline{a}',\underline{k}')$, hence $W_{p,q}(\underline{k})$
from $ W_{p,q'}(\underline{k}')$.

One verifies that the above criterion is independent of the choice
of the diffeomorphism $\phi$, whose ambiguities correspond to
Diff$^+(L(p,q))$. In fact, Lisca shows that even the
diffeomorphism type of the resulting manifold $V$ is independent
of the possible choices of $\phi$. For this he proves that
$W_{p,q}(\underline{k})$ admits an orientation-preserving
self-diffeomorphism which induces multiplication by $-1$ on
$H_1(L(p,q))$ (cf. \cite[(7.2)]{L 08}), and clearly
$\Pi(\underline{a})$  has a self-diffeomorphism which induces
multiplication by $q$ on $H_1(L(p,q))$ provided that
$q'=q$.

\subsection{Compatibility with the order}\label{ss:lisord}
One can eliminate the ambiguity left by (\ref{th:lis})(b) using
the order of the boundary. Notice that if $r=1$, or even if $r>1$
but both sequences $\underline{a}$ and $\underline{k}$ are
symmetric, then there is no ambiguity since $(p,q,\underline{k})=
(p,q',\underline{k}')$.

Assume that we are in the remaining situations. Recall that all
the time $\underline{a}$ and $q$ are related by the expansion
$[a_1,\ldots,a_r]=\frac{p}{p-q}$. The point is that the framed
link from Figure \ref{fig:Surgery2} is not symmetric.
If we mark the link components of ${\mathcal L}(\underline{k})$
as in Figure \ref{fig:Surgery2}, then by the rule described in subsection
(\ref{ss:4.3}) (namely, by imposing $\partial R_1\subset L_1$),
we appoint the preferred order of the boundary of $W_{p,q}(\underline{k})$,
 cf. (\ref{preferred}). The filling obtained in this way (the space
$W_{p,q}(\underline{k})$ with the preferred  order on its
boundary) will be  denoted by $W_{p,q}(\underline{k})^*$.

Notice that no orientation-preserving diffeomorphism
$W_{p,q}(\underline{k})^*\to
W_{p,q'}(\underline{k'})^*$ (from \linebreak (\ref{th:lis})(b))
preserves the preferred orders of the boundaries. Hence:

\begin{theorem}\label{th:lis*}
All the spaces  $W_{p,q}(\underline{k})^*$ are different, hence
their boundaries $L(p,q)^*$ and $\underline{k}\in
K_r(\underline{a})$ uniquely determine  all the Stein fillings up
to orientation-preserving diffeomorphisms which preserve the order
of the boundary.
\end{theorem}

The criterion (\ref{glueplumb}) will have the following new form.
Let $W^*$ be a Stein filling of $(L(p,q),\xi_{st})$ with an order
on its boundary. Consider $\Pi(\underline{a})^*$  with its
preferred order (providing a well-determined order of the
$s_i$'s). Construct $V$ as in (\ref{glueplumb}), and consider the
two pairs $(q,\underline{k})$ and $(q',\underline{k}')$ provided
(but undecided) by  (\ref{glueplumb}).

\begin{proposition}\label{glueplumb*} If $\phi$ preserves
  (resp. reverses) the orders of   the boundary then
    $W^*$ is orientation and order-preserving diffeomorphic to
    $W_{p,q}(\underline{k})^*$ (resp. to
     $W_{p,q'}(\underline{k}')^*$).
\end{proposition}

\begin{remark}
  Assume that $q=q'$. Then the permutation of the coordinates
  $(x,y)\mapsto (y,x)$ of  
$\C^2$ induces an automorphism of ${\mathcal X}_{p,q}$, and
also of the miniversal deformation space. The permutation on its
reduced components corresponds to  
$\underline{k}\mapsto \underline{k}'$ inducing an orientation
preserving diffeomorphism 
$W_{p,q}(\underline{k})\to W_{p,q}(\underline{k}')$. As it follows
from the above discussion, this diffeomorhism 
\emph{does not} preserve the order, provided that
$\underline{k}\not=\underline{k}'$. 
\end{remark}

\medskip
\section{The smoothings of cyclic quotient
  singularities, \\ following Christophersen and
  Stevens} \label{chrste}

In this section we recall some
 results of Christophersen and Stevens on the structure of the
reduced miniversal base space of cyclic quotients. For more
details see \cite{BR 95,C 91,S 03}.

\subsection{Generalities on versal deformations}\label{ss:6.1}

\begin{definition} \label{miniv}
  Let $(X,x)$ be a germ of a complex analytic space. A
  {\it deformation} of $(X,x)$ is a germ of flat morphism
  $\pi:(Y,y)\rightarrow (S,s)$ together with an
  isomorphism between $(X,x)$ and the special fiber $\pi^{-1}(s)$.
  A deformation of $(X,x)$ is {\it versal} if any other deformation is
  obtainable from it
  by a base-change. A versal deformation is {\it miniversal} if
  the Zariski tangent space of its base $(S,s)$ has the smallest possible
  dimension. A {\it smoothing component} is an irreducible component
  of the miniversal
  base space  over which the generic  fibers are smooth.
\end{definition}

If $(X,x)$ is a germ of a reduced complex analytic space
\emph{with an isolated singularity}, then the following well-known
facts hold:

  (i)  (Schlessinger \cite{S 68}, Grauert \cite{G 72}) The miniversal
  deformation $\pi$
  exists and is unique up   to (non-unique) isomorphism.

(ii) (Artin \cite{A 74}) If $(X,x)$ is a rational surface singularity,
then all the components of the reduced miniversal base space are
smoothing ones.  

(iii) (Looijenga \cite{L 84}) There exist (Milnor) representatives
$Y_{red}$ and $S_{red}$ 
of the reduced total and base spaces of $\pi$ such that the
restriction $\pi:\partial Y_{red} \cap\pi^{-1}(S_{red})
\rightarrow S_{red}$ is a trivial $C^{\infty}$-fibration.

\vspace{2mm}

Hence, for each smoothing component,  the oriented diffeomorphism type
of the oriented manifold with boundary
 $(\pi^{-1}(s)\cap Y_{red},\pi^{-1}(s)\cap \partial Y_{red})$ does not depend
 on the choice of the generic element   $s$: it
 is called \emph{the Milnor fiber} of that component. Moreover,
its boundary   is canonically
identified with the link \emph{up to isotopy}.  In particular, the Milnor
fiber over a
smoothing component is a Stein filling of the link endowed with
its standard contact structure (provided that the representatives
are carefully chosen, see \cite{CNP 06}).

\subsection{The equations of ${\mathcal X}_{p,q}$}\label{ss:6.2} Denote
by $S_{red}(p,q)$ the reduced base space of the miniversal
deformation of the cyclic quotient singularity
$\mathcal{X}_{p,q}$. Inspired by Arndt's work
 \cite{A 88}, Christophersen  defined in \cite{C 91} for
each $\underline{k}\in K_r(\underline{a})$ an explicit system
$\E_{\underline{k}}$ of equations which define
$\mathcal{X}_{p,q}$, and an explicit deformation
$\tilde{\mathcal{E}}_{\underline{k}}$ of these equations with
smooth parameter space. Based in an essential way on the work
\cite{KS 88} of Koll{\'a}r and Shepherd-Barron, Stevens proved in
\cite{S 91} that one gets in this way all the irreducible
components of $S_{red}(p,q)$. Hence:

\begin{theorem} \label{miniv2}
   The reduced base space $S_{red}(p,q)$ of the miniversal deformation of
   $\mathcal{X}_{p,q}$ has exactly $\# K_r(\underline{a})$ irreducible
   components. 
\end{theorem}

As $\mathcal{X}_{p,q}$ is a rational singularity, 
(\ref{ss:6.1})(ii) shows that all the irreducible components of
$S_{red}(p,q)$ are smoothing components. 
\emph{We denote by  $S^{CS}_{\underline{k}}$ the irreducible
component which corresponds to $\underline{k} \in
K_r(\underline{a})$ through the equations of Christophersen and Stevens.} 

\vspace{2mm}

Through the equations of Christophersen and Stevens one has in
fact an {\it explicit} bijection between the set
$K_r(\underline{a})$ and the irreducible components of
$S_{red}(p,q)$. This (together with (\ref{ss:5.2})) allows to
understand the meaning of Lisca's conjecture from the
introduction.

The system $\E_{\underline{k}}$  is best described using toric
  geometry.  The singularity $\mathcal{X}_{p, q}$ may also be seen as
  the germ at the   zero-dimensional orbit of the toric variety
  $\mathcal{Z}_{\sigma_{p,q}}= \mbox{Spec } \C
[\check{\sigma}_{p,q} \cap M]$, where $\sigma_{p,q}\subset
  N_{\R}$ is an oriented cone in $N$ of type $\frac{p}{q}$, and
$M:=\mbox{Hom}(N, \Z)$. We identify  $\check{\sigma}_{p,q}$ and
$\overline{\sigma}_{p,q}\simeq \sigma_{p, p-q}$, cf. Figure
\ref{fig:Suppl}. The lattice points
$(\overline{v}_0,\ldots,\overline{v}_{r+1})$ are the
  minimal generating set of the semigroup
  $\overline{\sigma}_{p,q} \cap N \: \simeq \: \check{\sigma}_{p,q} \cap M$.
  Therefore, the monomials:
\begin{equation} \label{defX}
  z_i:= \chi^{\overline{v}_i} \: \: \mbox{ for all } i \in \{0,\ldots,
  r+1\},
\end{equation}
generate the toric algebra $\C [\check{\sigma}_{p,q} \cap M]$.
Hence, the toric surface $\mathcal{Z}_{{\sigma}_{p,q}}$ may be
embedded inside $\C^{r+2}$ using the regular  functions
$z_0,\ldots,z_{r+1}$. Christophersen and Stevens write for each
$\underline{k} \in K_r(\underline{a})$ the special system
$\E_{\underline{k}}$ of binomial equations which define the image
$\mathcal{X}_{p,q}$ of $\mathcal{Z}_{\sigma_{p,q}}$ by this
embedding (i.e., the space $\mathcal{X}_{p,q}$ is independent of
$\underline{k}$, but for each $\underline{k}$ one  deforms
different set of equations of the space, this is the set
$\E_{\underline{k}}$). See \cite[pages 83-84]{C 91}, \cite[pages
316-317]{S 91} or \cite[pages 8-11]{BR 95} for different
presentations. Some of the equations are independent of
$\underline{k}$. They include:
\begin{equation} \label{simpleq}
z_{i-1}z_{i+1}-z_i^{a_i}=0\ \mbox{ for all } i \in \{1,\ldots,r\}.
\end{equation}
Using the special form of the equations, one defines their
deformations $\tilde{\mathcal{E}}_{\underline{k}}$, see \cite{C
91,S 91}.

\begin{remark} Once the preferred order of the coordinates
$(x,y)$ is fixed (cf. \ref{preferred} and \ref{markcoord}), they
also induce an order/marking of the coordinates $z_i$  via the
identities (\ref{defX}).
\end{remark}

\medskip
\section{The  smoothings of sandwiched surface
singularities, \\
following  de Jong and van Straten}
\label{topsmooth}

Cyclic quotient singularities are particular cases of
\emph{sandwiched surface singularities}. de Jong and van Straten
in \cite{JS 98} related the deformation theory of sandwiched
surface singularities to the deformation theory of so-called
\emph{decorated plane curve
  singularities}. They show that {\it 1-parameter deformations of
  decorated curves provide 1-parameter deformations for sandwiched
  singularities, and all of these later ones can be obtained in
  this way}. Moreover, the Milnor fibers of those  which are smoothings
  can be combinatorially described by  the so-called \emph{picture
    deformations}.

  In this section we explain the general framework, while in
  \S\,\ref{cycl}  we specialize it  to cyclic quotient singularities.

\subsection{Sandwiched singularities}\label{ss:7.1} The normal surface
singularity $(X,0)$ is called {\it sandwiched} if it is a germ of
an algebraic surface which admits a birational map $X\to \C^2$.
They  were introduced in \cite{S 90} by Spivakovsky; see also
\cite{JS 98,L 00,M 03} for different view-points.

 Sandwiched singularities are rational. They are
characterized (like the rational singularities) by their dual
resolution graphs. Hence one may  speak about \emph{sandwiched
graphs}.

\begin{proposition} \label{charsand}
Sandwiched graphs are characterized as follows: by adding new
vertices with weights (self-intersections ) $-1$ (on the `right
places')  one may obtain a `smooth graph', i.e.  the dual tree of
a configuration of ${\mathbb P}^1$'s which blows down to a smooth
point.
\end{proposition}

\subsection{Decorated curves and their deformations}\label{ss:7.2}
Any sandwiched singularity may be obtained from a weighted curve
$(C,l)$. Here $(C,0)\subset (\C^2,0)$ denotes a reduced germ of
plane curve with branches $\{C_i\}_{1 \leq i \leq r}$ and a
function $ \{1,\ldots,r\}\ni i\mapsto l_i\in \N^*$.

Consider the
minimal resolution of $C$. The \emph{multiplicity sequence}
associated with $C_i$ is the sequence of multiplicities on the
successive strict transforms of $C_i$, starting from $C_i$ itself
and not counting the last strict transform.  The {\it total
multiplicity} $m(i)$  of $C_i$ with respect to $C$ is the sum of
multiplicities of $C_i$ defined before.

\begin{definition} \cite[(1.3)]{JS 98} \label{decgerm}
  A {\it decorated germ} of plane curve is a weighted germ $(C,l)$
  such that $l_i\geq m(i)$ for all \ $  i \in \{1,\ldots,r\}$.
\end{definition}

The point is that starting from a decorated  germ, one can blow up
iteratively points infinitely near $0$ on the strict transform of
$C$, such that the number of such points sitting  on the strict
transform of $C_i$ is exactly $l_i$. If $l_i$ is sufficient large
(in general, larger than $m(i)$), then the union of the
exceptional components which do not meet the strict transform of
$C$ form a connected configuration of curves. After its
contraction one gets necessarily a sandwiched singularity
$X(C,l)$, determined uniquely by $(C,l)$ (for details see \cite{JS
98}).

The total multiplicity of $C_i$ with respect to $C$ may be encoded
also as the unique subscheme of length $m(i)$ supported on the
preimage of $0$ on the normalization of $C_i$. The same thing is
valid for $l_i$. This allows to define the \emph{total
multiplicity scheme} $m(C)$ of any reduced curve contained in a
smooth complex surface, as the union of the total multiplicity
schemes of all its germs:

\begin{definition} (i) \cite[(4.1)]{JS 98} \label{defcurve}
  Given a smooth complex analytic surface $\Sigma$, a  pair $(C,l)$ consisting
  of a reduced curve  $C\hookrightarrow \Sigma$ and a subscheme $l$ of
  the normalization $\tilde{C}$ of $C$ is called {\it a decorated
    curve} if $m(C)$ is a subscheme of \, $l$.

\vspace{1mm}

 (ii) \cite[p. 476]{JS 98}
  A $1$-{\it parameter deformation} of a decorated curve $(C,l)$
  over a germ of smooth curve $(S,0)$ consists of:

(1) a $\delta$-constant deformation $C_S \rightarrow S$ of $C$;

(2) a flat deformation $l_S \subset \tilde{C_S} =\tilde{C}\times
S$ of the scheme $l$, such that

(3) $m_S\subset l_S$, where the \emph{relative total multiplicity
  scheme} $m_S$ of $\tilde{C_S}\rightarrow C_S$ is defined as the
closure $\overline{\bigcup_{s \in S \setminus 0} m(C_s)}$.

\vspace{1mm}
 (iii)  A $1$-parameter deformation $(C_S, l_S)$ is called a
   {\it picture deformation} if for generic $s \neq 0$ the divisor
   $l_s$ is reduced.
\end{definition}

\cite[(4.4)]{JS 98} says that all the 1-parameter deformations of
$X(C,l)$ are obtained by $1$-parameter deformations of the
decorated germ $(C,l)$. Moreover, \emph{picture
  deformations} provide  smoothings of $X(C,l)$.

\subsection{Picture deformations and their Milnor fibers.}\label{ss:7.3}
Consider a decorated germ $(C,l)$ with {\it all the components
$\{C_i\}_{1\leq i\leq r} $ smooth}, and one of its picture
deformations $(C_S, l_S)$. Fix a closed Milnor ball $B$ for the
germ $(C,0)$. For $s\neq 0$ sufficiently small, $C_s$ will have a
representative in $B$, denoted by $D$, which meets $\partial B$
transversally. It is a union of embedded \emph{discs}
$\{D_i\}_{1\leq i\leq r}$ canonically oriented by their complex
structures (and whose set of indices correspond canonically to
those of  $\{C_i\}_{1\leq i\leq r} $). The singularities of $D$
consist of ordinary $m$-tuple points, for various $m$.

Denote by $\{P_j\}_{1 \leq j \leq n}$ the images in $B$ of the
points in the support of $l_s$. It is a finite set of points which
contains the singular set of $D$ (because $m_s \subset l_s$), but
it contains some other `free' points as well. There is a priori no
preferred choice of their ordering. Hence, the matrix introduced
next is well-defined only up to permutation of columns.

\begin{definition} \cite[page 483]{JS 98} \label{incid}
  The {\it incidence matrix} of a picture deformation $(C_S, l_S)$
  is the matrix $\mathcal{I}(C_S, l_S)\in \mbox{Mat}_{r,n}(\Z)$ whose
  entry at the  intersection
  of the $i$-th row and the $j$-th column is equal to $1$ if $P_j\in D_i$
  and to $0$ if $P_j \notin D_i$.
\end{definition}
\noindent Part (ii)(2) of (\ref{defcurve}) implies that:
\begin{equation} \label{nopoints}
  \mbox{the sum of entries on the $i$-th row of $\mathcal{I}(C_S, l_S)$ is
  $l_i$.}
\end{equation}
\noindent  The Milnor fiber of such a smoothing is recovered as
follows. Let:
\begin{equation} \label{blowpoints}
   (\tilde{B}, \tilde{D})
   \stackrel{\beta}{\longrightarrow}(B, D)
\end{equation}
  be the simultaneous blow-up of  the points $P_j$ of $D$. Here
  $\tilde{D}:= \cup_{1\leq i \leq r}\tilde{D}_i$, where $\tilde{D_i}$
  is the strict transform by the modification $\beta$ of the disc $D_i$.
Let $T_i$ be a sufficiently small open tubular neighborhood of
$\tilde{D}_i$ in $\tilde{B}$.

\begin{proposition}\cite[(5.1)]{JS 98} \label{recofiber}
  Suppose that all the irreducible components of $C$ are smooth. 
  The Milnor fiber of the smoothing of $X(C,l)$ corresponding to the
  picture deformation $(C_S, l_S)$ is orientation-preserving
  diffeomorphic to the
  compact oriented manifold with boundary $W:= \tilde{B}\setminus
         (\bigcup_{1\leq i \leq r} T_i)$ (whose corners are smoothed).
\end{proposition}

\medskip
\section{The smoothings of cyclic quotient singularities
\\ following de Jong and van Straten}
 \label{cycl}

\subsection{How to find $(C,l)$?}\label{ss:8.1} The construction of
the decorated germ $(C,l)$ used for cyclic quotients is in a
natural way valid for the more general class of {\it minimal
singularities}. Hence,  it is natural to present it in this
context.   Minimal singularities were introduced by Koll{\'a}r
\cite{K 85} in arbitrary dimension. In the case of normal
surfaces, they are exactly those rational germs which have reduced
fundamental cycle (see e.g. L{\^e} \cite{L 00}).  They are special
sandwiched singularities, characterized by their graphs as
follows.

Consider the minimal resolution of a rational singularity $(X,x)$.
Let $\Gamma$ be its dual graph, $J$ the set of its vertices,
$\nu_j$ the valency and $-e_j<0$ the weight (self-intersection) of
the vertex $j \in J$. Then $(X,x)$ is minimal if and only if its
minimal resolution satisfies $\nu_j\leq e_j$ for all $j\in J$.

For any minimal singularity $(X,x)$, there is an easy algorithm
which provides a decorated curve $(C,l)$ defining $(X,x)$,
starting from $\Gamma$. Its steps are the following:

\vspace{2mm}

I) \emph{Construct a new graph $\Gamma'$ by connecting the vertex
$j$ of $\Gamma$ to $e_j-\nu_j$ new vertices, with the exception of
one vertex $j_0$ satisfying $e_{j_0}-\nu_{j_0} >0$, which must be
connected to $e_{j_0}-\nu_{j_0}-1$ new vertices. Each new vertex
is endowed with the weight $-1$. }

II) \emph{Construct another graph $\Gamma''$ by endowing each new
vertex with one new arrowhead.}

\vspace{2mm}

Then $\Gamma'$ is a smooth graph  (hence by (\ref{charsand})
minimal singularities are sandwiched indeed), and  $\Gamma''$ is
the (not necessarily minimal) dual graph of an embedded resolution
of a germ of a plane curve $C$ with all components $C_i$ smooth.
They are obtained by blowing-down `curvettas' corresponding to the
arrowheads of $\Gamma''$.

The graph $\Gamma''$ has the following properties too:

\begin{lemma} \label{compl} 1) The intersection number of the irreducible
    curves $C_i$ and
    $C_j$ associated with two distinct arrowheads is equal to the number
    of vertices on the intersection of the geodesics from $j_0$ to the
    two arrowheads.

  2) The weight $l_i$ associated with the irreducible curve
  corresponding to an arrowhead is equal to the distance from the
  vertex $j_0$ to that arrowhead.
\end{lemma}

\begin{figure}[ht!]
\vspace*{6mm}
\labellist
\small\hair 2pt
\pinlabel $j_0$ at 5 110
\pinlabel  {$m_1 -1$} at 32 -20
\pinlabel $m_2$ at 210 -20
\pinlabel $m_t$ at 500 -20
\pinlabel $m_{t+1}$ at 681 -20
\pinlabel  {$\dots$} at 34 40
\pinlabel $\dots$ at 212 40
\pinlabel $\dots$ at 502 40
\pinlabel $\dots$ at 683 40
\pinlabel {$n_1 -3$} at 123 52
\pinlabel {$n_t -3$} at 593 52
\pinlabel $-2$ at 92 126
\pinlabel $-2$ at 151 126
\pinlabel $-2$ at 271 126
\pinlabel $-2$ at 439 126
\pinlabel $-2$ at 558 126
\pinlabel $-2$ at 618 126
\pinlabel {$-(m_1 +1)$} at 33 130
\pinlabel {$-(m_{2}+2)$} at 212 130
\pinlabel {$-(m_{t}+2)$} at 501 130
\pinlabel {$-(m_{t+1}+1)$} at 688 130
\pinlabel {$C_1$} at 0 32
\pinlabel {$C_{m_1-1}$} at 82 32
\pinlabel {$C_{m_1}$} at 172 32
\pinlabel {$C_{m_1 + m_2 -1}$} at 275 32
\pinlabel $-1$ at 5 74
\pinlabel $-1$ at 65 74
\pinlabel $-1$ at 185 74
\pinlabel $-1$ at 245 74
\pinlabel $-1$ at 465 74
\pinlabel $-1$ at 530 74
\pinlabel $-1$ at 650 74
\pinlabel $-1$ at 709 74
\endlabellist
\centering
\includegraphics[scale=0.60]{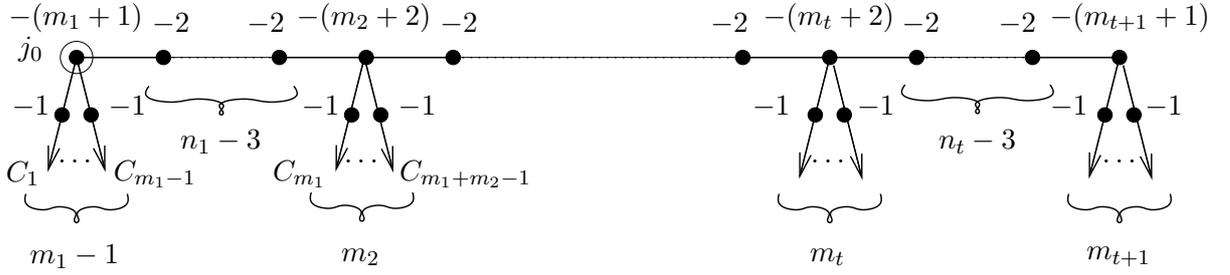}
\vspace{2mm}
\caption{The dual graph associated with $(C,l)$ with
  $\mathcal{X}_{p,q}=X(C,l)$}
\label{fig:Diacurve}
\end{figure}

We apply the previous procedure to $\mathcal{X}_{p,q}$, by
choosing as $j_0$ that vertex of $\Gamma=G(\underline{b})$ which
is marked by 1, i.e.  whose weight is $-b_1$ (see Figure
\ref{fig:Gb}). We use the notations  (\ref{HJexp2}) for the
sequence $\underline{b}$.  Hence, $\mathcal{X}_{p,q}$ may be
presented as a sandwiched singularity $X(C,l)$, where the
components of $C$ are indicated in the graph $\Gamma''$ shown in
Figure \ref{fig:Diacurve}. By (\ref{eq:S4}), the number of curves
is exactly $r$ (which explains why we have chosen this notation for
the number of components of $C$ in the previous section). By
(\ref{compl}) one has:
 \begin{equation} \label{inters}
     C_i\cdot C_j=l_i -1 \ \ \ \mbox{for all} \ \ \   1 \leq i < j \leq
     r;
 \end{equation}
\begin{equation} \label{li}
   l_i=2+\sum_{1\leq j\leq h-1} \, (n_j-2)  \ \ \ \mbox{ whenever } \ \ \
   \sum_{1\leq j\leq h-1}m_j\leq i < \sum_{1\leq j\leq h}m_j.
 \end{equation}
By (\ref{HJexp1}), these last relations (\ref{li}) transform into:
\begin{equation} \label{la}
   l_i=2+ \sum_{1\leq j\leq i} (a_i-2) \ \ \mbox{for all} \ \ i \in
   \{1,\ldots,r\}.
\end{equation}

\vspace*{4mm}

\subsection{From triangulations to the incidence matrix}\label{ss:8.2}
We will describe the {\bf incidence matrix} 
using an interpretation by Stevens of the elements of $K_r$ via
triangulations of a polygon \cite{S 91}, and the following
notation: If $M\in\mbox{Mat}_{r,n}(\Z)$ then $\smallint
M\in\mbox{Mat}_{r,n}(\Z)$ will denote that matrix whose $i$-th row
is the sum of the first $i$ rows of $M$. By our construction, we
get the {\it same} incidence matrix as in \cite{JS 98},
nevertheless we arrive to it slightly differently (perhaps, more
conceptually). \cite{JS 98} starts with a `difference matrix'
(with non-negative entries), considers its $\smallint$, and the
$modulo \ 2$ remainders of the entries of this second matrix
constitute the incidence matrix. In our case, we conceptually
assign to each entry of the `difference matrix' a sign such that
its $\smallint$ will be exactly the incidence matrix; cf.
(\ref{re:JSmod2}).

Assume $r>1$ (since $\#K_1=1$, in the identifications we wish to
get we lose nothing).

Consider a convex polygon $\mathcal{P}_{r+1}$ in the plane with
$r+1$ vertices marked successively by $A_1,\ldots,A_{r+1}$. Denote
by $T(\mathcal{P}_{r+1})$ the set of triangulations of
$\mathcal{P}_{r+1}$ with vertices $A_1,\ldots, A_{r+1}$ and edges
which are diagonals of $\mathcal{P}_{r+1}$. With each
triangulation $\theta\in T(\mathcal{P}_{r+1})$ associate the
sequence $(k_1,\ldots,k_r)$ such that $k_i$ is the number of
triangles containing the vertex $A_i$. Then $\theta\mapsto
\underline{k}$ realizes a bijection between $T(\mathcal{P}_{r+1})$
and $K_r$, cf. \cite{S 91}.

Fix $\theta\in T(\mathcal{P}_{r+1})$. To each triangle $\Delta$ of
$\theta$, and vertex $A$ of  $\mathcal{P}_{r+1}$ we define a
`sign' $\alpha =\alpha(A,\Delta)\in\{0,-1,+1\}$ as follows. If $A$
is not a vertex of $\Delta$, we take $\alpha=0$. Then  we
(totally) order the vertices of $\Delta$ by restricting  to them
the order $A_1,\ldots, A_{r+1}$ of the vertices of
$\mathcal{P}_{r+1}$. We set $\alpha=+1$ for the first and third
vertex, while $\alpha=-1$ for the second one. Next, order the
triangles $\{\Delta_j\}_{j=1}^{r-1}$ of $\theta$ in an arbitrary
way, and define the   `sign-incidence matrix' between the vertices
$A_1,\ldots, A_r$ (corresponding to the rows) and the triangles
$\Delta_1,\ldots, \Delta_{r-1}$ (corresponding to the columns) by
$d_{i,j}: =\alpha(A_i,\Delta_j)$. If $\theta$ corresponds to
$\underline{k}$, then denote this matrix by $D(\underline{k})\in
\mbox{Mat}_{r, r-1}( \Z)$ (well-defined up to a permutation of its
columns). Moreover, for each $\ell\in \N$, denote by
$M_{r,\ell}(i)\in \mbox{Mat}_{r, \ell}( \Z)$
 the matrix whose
 entries are all equal to 1 on the $i$-th row and all equal to $0$
 elsewhere.

Then,  for each $\underline{k}\in K_r(\underline{a})$, consider
the `block-matrix':
\begin{equation} \label{bloc}
  D(\underline{a}; \underline{k}):=
    (\: D(\underline{k}) \: | \: M_{r,a_1-k_1}(1) \:  | \:
    \cdots \: | \:  M_{r, a_r-k_r}(r) \:) \: \in \:
     \mbox{Mat}_{r, r-1+\sum( a_i -  k_i)}( \Z).
\end{equation}
The following theorem, valid for any fixed $(C,l)$ as in
(\ref{ss:8.1}), follows from \cite[Section 6.4]{JS 98}:

  \begin{theorem}  \cite[(6.18)]{JS 98} \label{expict}
    For every $\underline{k}\in     K_r(\underline{a})$, the matrix
      $\smallint D(\underline{a}; \underline{k})$ is (up to a permutation
      of columns) the incidence matrix
      of some picture deformation of $(C,l)$.
In particular, the
    number $n$  of points $\{P_i\}_i$ is
$ n=r-1+\sum_{i=1}^r (a_i-k_i)$. Moreover,
      varying $\underline{k}\in  K_r(\underline{a})$, one gets all
      incidence matrices
    of picture     deformations of $(C,l)$.
  \end{theorem}

 \begin{remark}\label{re:JSmod2}
    In \cite{JS 98} the authors call the
    matrices  $\smallint D(\underline{a}; \underline{k})$
    \emph{CQS-matrices}, and denote them by $M$. Up to the signs of the
    entries,  $D(\underline{a}; \underline{k})$ are  their
    \emph{difference matrices}
    $\Delta M$. More precisely, the entries of $\Delta M$
    are the absolute values
    of the entries of $D(\underline{a}; \underline{k})$.
 \end{remark}

\begin{remark}
  By (\ref{incid}) and (\ref{nopoints}), an incidence
  matrix has all its entries equal to $0$ or $+1$ and the sum of
  all entries of the $i$-th row is $l_i$. Let us verify that
  this is indeed the case for  $\smallint D(\underline{a};
  \underline{k})$. The first property is a consequence of the fact
  that on each column of $D( \underline{k})$, the
  non-zero entries are either $(+1, -1)$ or $(+1, -1, +1)$, always
  in this order, depending on the fact that $A_{r+1}$ is a vertex of the
  corresponding triangle or not. The second property is a consequence
  of (\ref{la}) and the following elementary fact regarding the
  above sign-assigning procedure: $A_1$ has no sign equal to $-1$,
  and $A_i$ has exactly one sign equal to $-1$ for any $i>1$.
\end{remark}

Since via picture deformations we hit all the  components of the
reduced miniversal base space of $\mathcal{X}_{p,q}$, via the
correspondence (\ref{expict}) de Jong and van Straten parametrize
the smoothing components by  the elements of $K_r(\underline{a})$.
Denote by $S_{\underline{k}}^{JS}$ {\it the component parametrized
by} $\underline{k} \in K_r(\underline{a})$.  Denote by  $
W(\underline{a}, \underline{k})$ {\it the Milnor fiber} of the
corresponding smoothing, i.e. the manifold constructed in
(\ref{ss:7.3}), specialized to the present situation.

\begin{example}
Here we list all the  objects presented in this section, applied
for  the singularity $\mathcal{X}_{11,4}$, i.e. for
$\underline{a}=(2,3,2,2)$ and $\underline{b}=(3,4)$. The dual
graph of the minimal resolution is $G(3,4)$, while the resolution
graph of the decorated germ $(C,l)$ is the following (where in
parenthesis we inserted the integers $\{l_i\}_{i=1}^4$) :

\vspace{3mm}
\begin{figure}[h!]
\labellist
\small\hair 2pt
\pinlabel  $-3$ at 5 131
\pinlabel  $-4$ at 220 131
\pinlabel  $-1$ at 28 37
\pinlabel  $-1$ at 142 37
\pinlabel  $-1$ at 200 37
\pinlabel  $-1$ at 290 37
\pinlabel  $(2)$ at 5 -20
\pinlabel  $(3)$ at 148 -20
\pinlabel  $(3)$ at 221 -20
\pinlabel  $(3)$ at 291 -20
\endlabellist
\centering
\includegraphics[scale=0.40]{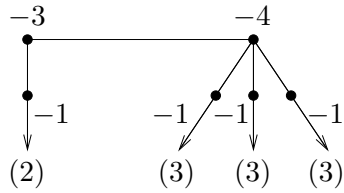}
\vspace{2mm} \caption{Dual graph associated with $(C,l)$, with
$\mathcal{X}_{11,4}=X(C,l)$} \label{fig:Exdualbis}
\end{figure}

\noindent
The associated set of sequences representing zero is:
$$K_5(2,3,2,2)=\{(1,2,2,1), (1,3,1,2)\}.$$

One has the following associated triangulations of a pentagon
(with the corresponding signs),  matrices
$D(\underline{a};\underline{k}), \smallint
D(\underline{a};\underline{k}) $ (with blocks $M_{r, a_i-k_i}$
separated, those with $a_i=k_i$ being empty), and generic member
of a picture deformation with $\smallint
D(\underline{a};\underline{k}) $ as incidence matrix :

\medskip
$\bullet$ $\underline{k}=(1,2,2,1)$.

\begin{figure}[h!]
\labellist
\small\hair 2pt
\pinlabel  $A_1$ at 237 75
\pinlabel  $A_2$ at 162 210
\pinlabel  $A_3$ at  55 210
\pinlabel  $A_4$ at  -30 75
\pinlabel  $+1$ at  188 76
\pinlabel  $-1$ at 166 135
\pinlabel  $+1$ at 127 162
\pinlabel  $-1$ at 81 164
\pinlabel  $+1$ at 48 132
\pinlabel  $-1$ at  23 76
\pinlabel  $\Delta_1$ at  145 52
\pinlabel  $\Delta_2$ at  108 89
\pinlabel  $\Delta_3$ at  70 52
\endlabellist
\centering
\includegraphics[scale=0.45]{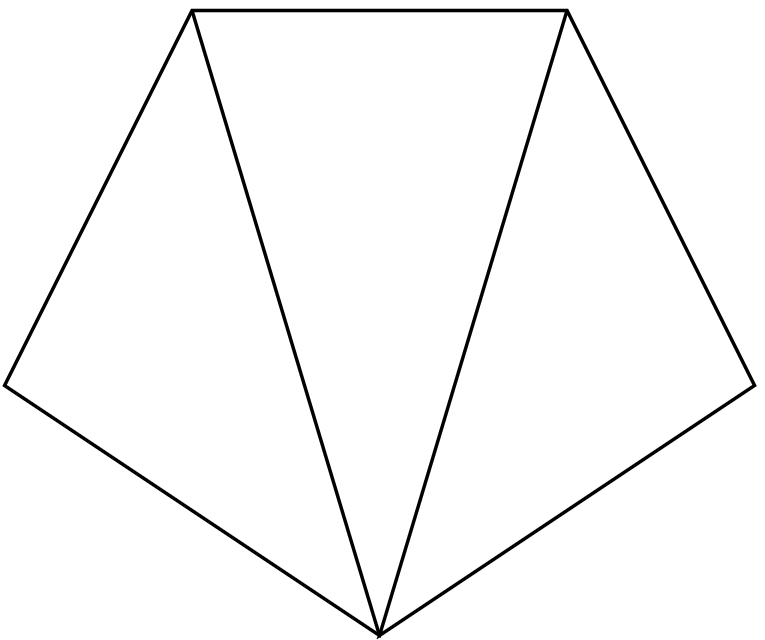}
\label{fig:Polyg1}
\end{figure}

$$ D(2,3,2,2 \: ; \: 1,2,2,1)=
  \left(  \begin{array}{ccc}
              +1 & 0 & 0 \\
              -1 & +1 & 0 \\
              0 & -1 & +1 \\
              0 & 0 & -1 \\
           \end{array} \right.
  \left| \begin{array}{c}
              +1 \\ 0 \\ 0 \\ 0
           \end{array} \right.
  \left| \begin{array}{c}
              0 \\ +1 \\ 0 \\ 0
           \end{array} \right|
 \left| \begin{array}{c}
              0 \\ 0 \\ 0 \\ +1
           \end{array} \right)
$$
$$ \smallint D(2,3,2,2 \: ; \: 1,2,2,1)=
  \left(  \begin{array}{cccccc}
              1 & 0 & 0 & 1 & 0 & 0 \\
              0 & 1 & 0 & 1 & 1 & 0 \\
              0 & 0 & 1 & 1 & 1 & 0 \\
              0 & 0 & 0 & 1 & 1 & 1\\
           \end{array}  \right)
$$

%\vspace{5mm}
\begin{figure}[h!]
\labellist
\small\hair 2pt
\pinlabel  $D_1$ at 259 57
\pinlabel  $D_2$ at 225 252
\pinlabel  $D_3$ at  137 276
\pinlabel  $D_4$ at  23 244
\pinlabel  $P_1$ at  61 79
\pinlabel  $P_4$ at 155 74
\pinlabel  $P_2$ at 67 130
\pinlabel  $P_3$ at 112 151
\pinlabel  $P_6$ at 186 130
\pinlabel  $P_5$ at  150 197
\endlabellist
\centering
\includegraphics[scale=0.35]{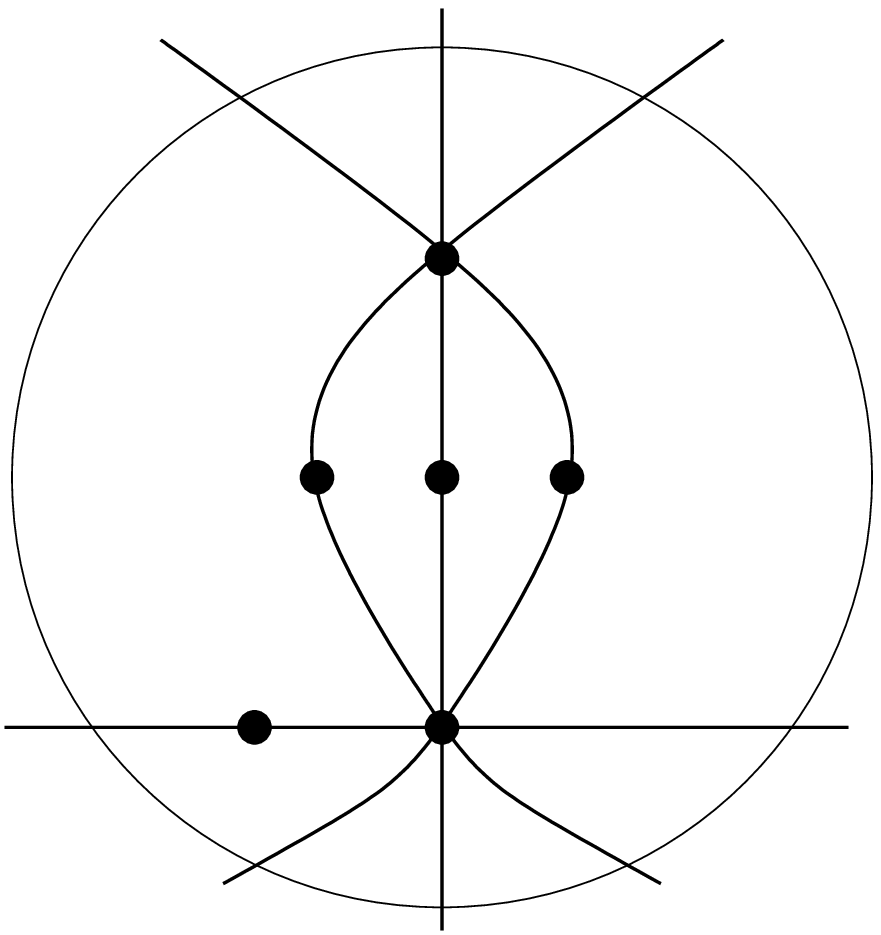}
%\vspace{2mm}
\caption{A general member of a picture deformation with
  $\underline{k}=(1,2,2,1)$}
\label{fig:Picture1}
\end{figure}

\medskip

%\pagebreak
$\bullet$ $\underline{k}=(1,3,1,2)$.

%\vspace{3mm}
\begin{figure}[h!]
\labellist
\small\hair 2pt
\pinlabel  $A_1$ at 237 75
\pinlabel  $A_2$ at 162 210
\pinlabel  $A_3$ at  55 210
\pinlabel  $A_4$ at  -30 75
\pinlabel  $+1$ at  190 82
\pinlabel  $-1$ at 166 135
\pinlabel  $+1$ at 135 150
\pinlabel  $-1$ at 37 82
\pinlabel  $+1$ at 123 170
\pinlabel  $-1$ at  73 170
\pinlabel $+1$ at 37 116
\pinlabel  $\Delta_1$ at  150 62
\pinlabel  $\Delta_2$ at  98 92
\pinlabel  $\Delta_3$ at  67 145
\endlabellist
\centering
\includegraphics[scale=0.45]{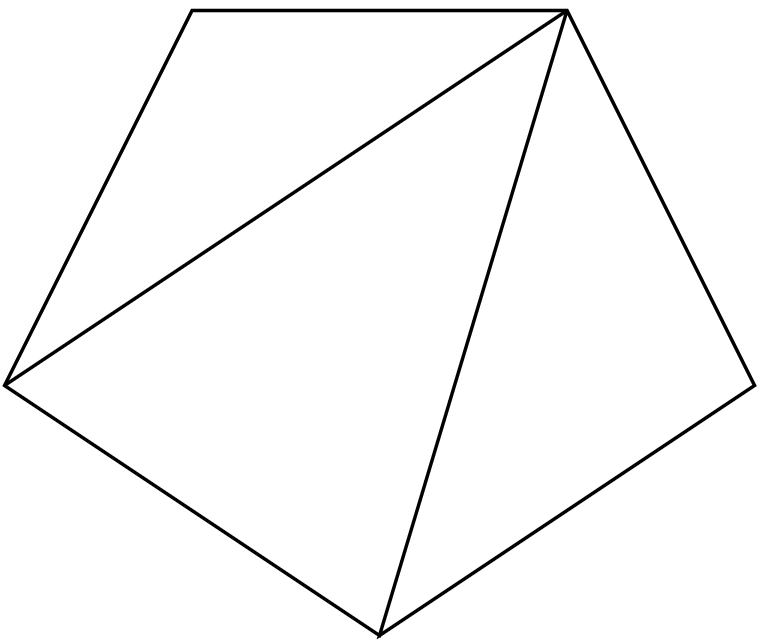}
\label{fig:Polyg2}
\end{figure}

$$ D(2,3,2,2 \: ; \: 1,3,1,2)=
  \left(  \begin{array}{ccc}
              +1 & 0 & 0 \\
              -1 & +1 & +1 \\
              0 & 0 & -1 \\
              0 & -1 & +1 \\
           \end{array} \right.
  \left| \begin{array}{c}
              +1 \\ 0 \\ 0 \\ 0
           \end{array} \right|
  \left| \begin{array}{c}
              0 \\ 0 \\ +1 \\ 0
           \end{array} \right|
   \left| \begin{array}{c}
               \\  \\  \\ \\
           \end{array} \right)
$$

$$ \smallint D(2,3,2,2 \: ; \: 1,3,1,2)=
  \left(  \begin{array}{ccccc}
              1 & 0 & 0 & 1 & 0 \\
              0 & 1 & 1 & 1 & 0 \\
              0 & 1 & 0 & 1 & 1 \\
              0 & 0 & 1 & 1 & 1 \\
           \end{array}  \right)
$$

\vspace{5mm}
\begin{figure}[h!]
\labellist
\small\hair 2pt
\pinlabel  $D_1$ at 259 57
\pinlabel  $D_2$ at 225 252
\pinlabel  $D_3$ at  137 276
\pinlabel  $D_4$ at  7 230
\pinlabel  $P_1$ at  61 79
\pinlabel  $P_4$ at 155 74
\pinlabel  $P_2$ at 148 194
\pinlabel  $P_3$ at 71 161
\pinlabel  $P_5$ at  152 156
\endlabellist
\centering
\includegraphics[scale=0.35]{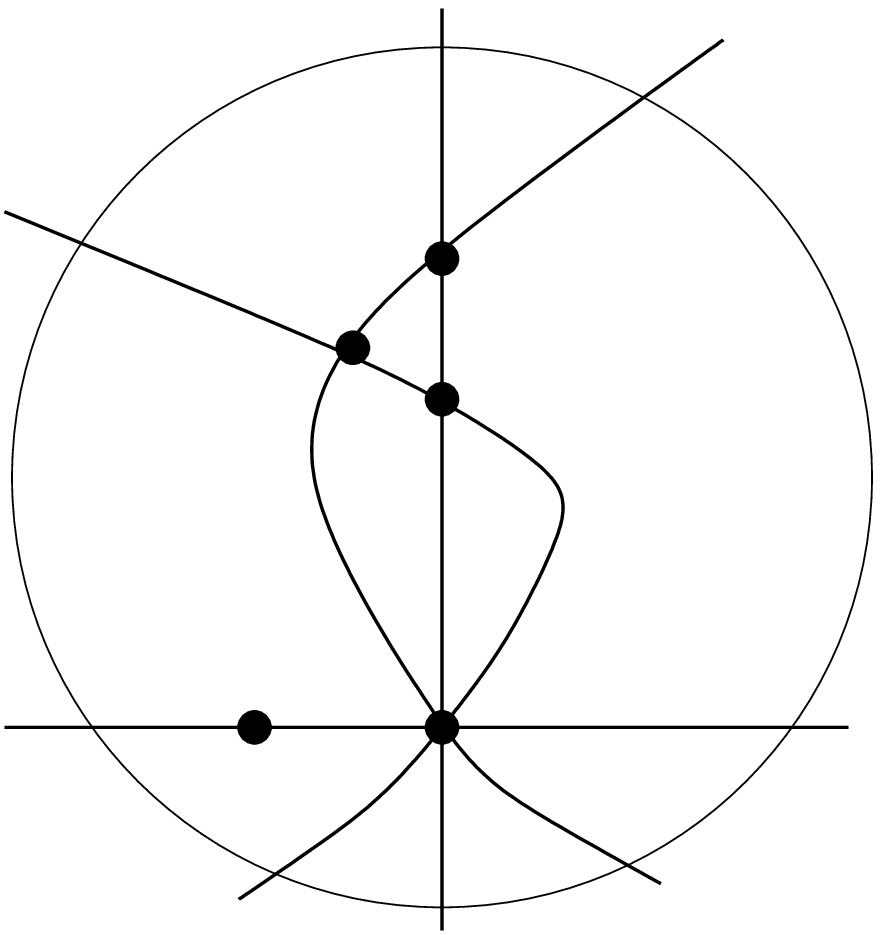}
%\vspace{2mm}
\caption{A general member of a picture deformation with
  $\underline{k}=(1,3,1,2)$}
\label{fig:Picture2}
\end{figure}
\end{example}

\section{`From Christophersen and Stevens to Lisca'}
\label{lcs}

\subsection{Inequalities} Before we start our discussion regarding
deformations of cyclic quotients, we state a technical lemma,
which will be used several times.

\begin{lemma}\label{lem:CC'} Assume that $\underline{x}\in \N^n$
is an admissible sequence (cf. (\ref{admis})).
 Then:

(1) For any $x_i'\geq x_i$ ($1\leq i\leq n$), $\underline{x}'$ is
admissible too.
% and $Z_n(x_1',\ldots,x_n')\geq  Z_n(x_1,\ldots,x_n)$.

(2) Assume that $\{\nu_i\}_{i=0}^{n+1}$ satisfy the inequalities
$\nu_{i+1}\geq x_i\, \nu_i-\nu_{i-1}$ for all $1\leq i\leq n$. Then for
all $i$ one also has:
\begin{equation}\label{eq:C3}
\nu_{i+1}\geq Z_i(x_i,\ldots,x_1)\,\nu_1-Z_{i-1}(x_i,\ldots,
x_2)\,\nu_0.
\end{equation}

(3) Assume that  $x_i\geq 2$ for all $i$ and set
$\mu_i:=Z_{i-1}(x_1,\ldots, x_{i-1})-Z_{i-2}(x_2,\ldots,x_{i-1})$
for $1\leq i\leq n$. Then
\begin{equation*}
1=Z_0<Z_1(x_1)<Z_2(x_1,x_2)<\cdots < Z_n(x_1,\ldots,x_n)\,, \ \mbox{and}
\end{equation*}%and
\begin{equation*}
1=\mu_1\leq \mu_2\leq \cdots\leq \mu_n\,.
\end{equation*}
\end{lemma}

\begin{proof} (1) 
$M(\underline{x}')$ is the sum of
$M(\underline{x})$ and a diagonal matrix with non-negative
entries. Therefore, if $M(\underline{x})$ is positive semi-definite of
rank $\geq (n-1)$, then $M(\underline{x}')$ will have the same
property.

(2)  We prove by decreasing induction on $j$ that for fixed $i$
and for all $ 1\leq j\leq i$ one has
\begin{equation}\label{eq:C2}
\nu_{i+1}\geq Z_{i-j+1}(x_i,x_{i-1},\ldots,
x_j)\,\nu_j-Z_{i-j}(x_i,x_{i-1},\ldots x_{j+1})\,\nu_{j-1}.
\end{equation}
For $j=i$ this is clear. The induction runs as
follows. Since $M(\underline{x})$ is positive semi-definite
 $Z_{i-j+1}(x_i,\ldots,x_j)\geq 0$, hence the right
hand side of (\ref{eq:C2}) is greater than:
$$
Z_{i-j+1}(x_i,\ldots,x_j)\,(x_{j-1}\nu_{j-1}-\nu_{j-2})-Z_{i-j}(x_i,\ldots
x_{j+1})\,\nu_{j-1}=
$$
$$[x_{j-1}Z_{i-j+1}(x_i,\ldots,x_j)-Z_{i-j}(x_i,\ldots
x_{j+1})]\,\nu_{j-1}-Z_{i-j+1}(x_i,\ldots,x_j)\,\nu_{j-2}.
$$
Since in the parenthesis we have exactly
$Z_{i-j+2}(x_i,\ldots,x_{j-1})$, cf. (\ref{recpol}), 
the proof of (\ref{eq:C2}) is
finished for all $j$. For $j=1$ we get the wished inequality.

(3) The first part follows by induction: if $Z_{i-1}(x_1,\ldots,
x_{i-1})>Z_{i-2}(x_1,\ldots, x_{i-2})$, then
$$Z_i(x_1,\ldots, x_i)=x_iZ_{i-1}(x_1,\ldots,
x_{i-1})-Z_{i-2}(x_1,\ldots, x_{i-2})>(x_i-1)Z_{i-1}(x_1,\ldots,
x_{i-1}).$$ This reinterpreted also shows that $\mu_i\geq 0$. Then
the identity $\mu_{i+1}-\mu_i=(x_i-2)\mu_i+(\mu_i-\mu_{i-1})$ and
induction finishes the proof.
\end{proof}
Since $\underline{k}$ is admissible and $\underline{k}\leq
\underline{a}$, the above lemma can be applied for both
$\underline{k}$ and $\underline{a}$.

\subsection{$\mathcal{X}_{\underline{k}}^t$ as the Milnor
  fiber.}\label{ss:9.1}  In
the sequel we will  follow the notations of \S \ref{chrste}. First
we concentrate on  $\xpq$.  Using equations (\ref{simpleq}) and
induction, one shows  that the restriction of each $z_i$  to
$\xpq$ is a rational function  in $(z_0,z_1)$ of the form:
\begin{equation}\label{eq:9.1}
z_i=z_1^{Z_{i-1}(a_1,\ldots, a_{i-1})}\cdot
z_0^{-Z_{i-2}(a_2,\ldots,a_{i-1})} \ \ \mbox{for} \ i\in
\{1,\ldots, r+1\}.
\end{equation}

The equations $\mathcal{E}_{\underline{k}}$ are weighted
homogeneous, however the weights $w_i:=w(z_i)$ are not unique.
With the choice  $w_0=w_1=1$ one has:
\begin{lemma}\label{lem:9.1}
(a) \ $ w_i=Z_{i-1}(a_1,\ldots, a_{i-1})-Z_{i-2}(a_2,\ldots,a_{i-1})$ for all
$i\geq 1$,

(b) \  $ 1=w_0=w_1\leq w_2\leq \cdots \leq w_{r+1}=q$.
\end{lemma}

\begin{proof}
(a) follows from (\ref{eq:9.1}), $w_{r+1}=q$ from
(\ref{contquot}), the rest of (b) from (\ref{lem:CC'}).
\end{proof}

We consider a special 1-parameter deformation
$\E_{\underline{k}}^t$ of the equations $\E_{\underline{k}}$. This
deformation is uniquely determined by the deformed equations of
(\ref{simpleq}) (cf. \cite{C 91}, \cite[(2.2)]{S 91}).  These last are:
\begin{equation} \label{eq:9.2}
z_{i-1}z_{i+1}=z_i^{a_i}+t\cdot z_i^{k_i}\ \mbox{ for all } i \in
\{1,\ldots,r\},
\end{equation}
where $t\in \C$. Note that, although (\ref{simpleq}) did not depend
on $\underline{k}$, this is not the case for their
deformations. Let $\mathcal{X}_{\underline{k}}^t$ be the affine space
determined
by the equations $\E_{\underline{k}}^t$ in $\C^{r+2}$.

\begin{lemma}\label{lem:9.2} The deformation $t\mapsto
  \mathcal{X}_{\underline{k}}^t$ has negative weight and
 is a smoothing belonging to the component
$S_{\underline{k}}^{CS}$. In particular,
$\mathcal{X}_{\underline{k}}^t$ is a smooth affine variety for
$t\not=0$.
\end{lemma}
\begin{proof}
The first statement just means that the weight of the added
monomial $z_i^{k_i}$ is not larger than the weight of
$z_{i-1}z_{i+1}-z_i^{a_i}$, i.e. $w_i>0$ and $k_i\leq a_i$. For
the second part one checks the general form of the equations of
$S_{\underline{k}}^{CS}$ from \cite[(2.2)]{S 91} or \cite{C 91},
and the fact that the present deformation does not belong to the
discriminant of $S_{\underline{k}}^{CS}$ described in  \cite{C
91}.

In fact, the smoothness  also follows from our direct computation,
as a byproduct of (\ref{indetfin}). Indeed,
$\mathcal{X}_{\underline{k}}^t$ as a fiber of the miniversal
deformation is normal. In (\ref{ss:9.3}) we will construct a
resolution of it, which has no exceptional curve by
(\ref{indetfin}). Hence $\mathcal{X}_{\underline{k}}^t$ is smooth.
\end{proof}

In particular, the above smoothing has a series of pleasant
properties. (E.g., it induces a projective deformation which is
locally trivial `near $\infty$'.) Moreover, by \cite[(2.2)]{W 81}:
\begin{equation}\label{eq:9.3}
\mathcal{X}_{\underline{k}}^t \ \mbox{is diffeomorphic to the Milnor fiber of
  $S_{\underline{k}}^{CS}$.}
\end{equation}

In the sequel we will denote by
$\widehat{\mathcal{X}_{\underline{k}}^t}$ the closure of
$\mathcal{X}_{\underline{k}}^t$ in $\PP^{r+2}$, and
let
$C_{\underline{k}}^\infty=
%\widehat{\mathcal{X}_{\underline{k}}^t}\cap H^\infty=
\widehat{\mathcal{X}_{\underline{k}}^t}\setminus
\mathcal{X}_{\underline{k}}^t$ be its     curve at infinity (as
may be seen by a computation, or by the equisingularity at
infinity mentioned above, $C_{\underline{k}}^\infty$ is
topologically independent of $\underline{k}$ and $t$).

\subsection{$\mathcal{X}_{\underline{k}}^t$ as a rational
  surface.}\label{ss:9.3}
Similarly as for $\xpq$ one shows by induction that on
$\mathcal{X}_{\underline{k}}^t$ all the restrictions of the
coordinates $z_i$ can be expressed as rational functions in
$(z_0,z_1)$:

\begin{lemma}\label{lem:9.3}
  For each $i \in \{1,\ldots,r+1\}$,  on
  $\mathcal{X}_{\underline{k}}^t$ one has:
 \begin{equation}\label{eq:9.5}
 z_i =z_0^{-Z_{i-2}(a_2,\ldots,a_{i-1})} P_i
\end{equation}
for some $P_i \in \Z [t,z_0, z_1]$. The polynomials $P_i$ satisfy
the inductive relations:
\begin{equation}\label{eq:9.6}
P_{i-1}\cdot P_{i+1}=P_i^{a_i}+tP_i^{k_i}\cdot z_0^{(a_i-k_i)\cdot
Z_{i-2}(a_2,\ldots,a_{i-1} )}
\end{equation}
with $P_1=z_1$ and with the convention $P_0=1$. Moreover $z_0\nmid
P_i$.
\end{lemma}

\begin{proof} Define $P_i$ by (\ref{eq:9.5}).
By a substitution it is clear that (\ref{eq:9.6}) follows from
(\ref{eq:9.2}) and (\ref{eq:9.5}). By (\ref{eq:9.6}) and
induction, $P_i$ is a (a priori rational) function in the variables
$(t,z_0,z_1)$. Hence, we only have to prove that $P_i$ is a
polynomial and $z_0\nmid P_i$. Let $R$ be an irreducible
polynomial in $(t,z_0,z_1)$, and $\nu_R:\C(t,z_0,z_1)^*\to\Z$ be
the valuation associated with it. We have to show that
$\nu_R(z_i)\geq 0$ for $R\not=z_0$ and $\nu_{z_0}(z_i)=
-Z_{i-2}(a_2,\ldots, a_{i-1})$.

Set $\nu_i:=\nu_R(z_i)$, and consider first $R=z_0$.
Then analysing (\ref{eq:9.2}) we get that 
$z_0$ is a pole of $z_i$ for $i\geq 2$,
 hence $\nu_{z_0}(z_i^{a_i}+tz_i^{k_i})\geq 
\nu_{z_0}(z_i^{a_i})$. This shows that  
$\nu_{i+1}\geq a_i\, \nu_i-\nu_{i-1}$, hence (\ref{eq:C3}) can be
applied. 
Since   $\nu_0=1$ and $\nu_1=0$ we get 
$\nu_i\geq -Z_{i-2}(a_2,\ldots, a_{i-1})$.

If $R\not=z_0$  then $\nu_0=0$
and $\nu_1\geq 0$. Assume that $\nu_j\geq 0$ for 
$0\leq j\leq i$. Then by (\ref{eq:9.2}) 
$\nu_{j+1}\geq k_j\, \nu_j-\nu_{j-1}$ for all $1\leq j\leq i$, hence
 (\ref{eq:C3}) can again be
applied (whose `$Z$-coefficients' are non-negative by the
admissibility of $\underline{k}$). In particular $\nu_{i+1}\geq 0$ too. 
Hence $P_i$ is a polynomial. Finally, (\ref{eq:9.6}) shows that
$P_{i+1}P_{i-1}\equiv cP_i^{a_i}$ (mod $z_0$) for some non-zero
constant $c$. Since $z_0$ does not divide $P_0$ and $P_1$, by
induction it does not divide $P_i$ either.
\end{proof}
In fact, by  (\ref{lem:CC'})(3) and (\ref{contquot}),  the
different exponents of $z_0$ in (\ref{eq:9.5}) satisfy
\begin{equation}\label{eq:EGY}
1=Z_0<Z_1(a_2)<\cdots < Z_{i-2}(a_2,\ldots, a_{i-1})<\cdots<
Z_{r-1}(a_2,\ldots, a_r)=p-q.
\end{equation}
Define now the application $\pi: \C^2\setminus
\{z_0=0\}\longrightarrow \mathcal{X}_{\underline{k}}^t$ \,  by
$(z_0,z_1)\mapsto (z_0,z_1,\ldots,z_{r+1})$, or
\begin{equation}\label{eq:9.7}
(z_0,z_1)\longrightarrow (z_0,z_1,z_0^{-1} P_2,\ldots,
z_0^{-Z_{i-2}(a_2,\ldots,a_{i-1})} P_i, \ldots,
z_0^{-(p-q)}P_{r+1})\in\C^{r+2},\end{equation} and  the induced
birational map 
$\widehat{\pi}:\PP^2\dashrightarrow
\widehat{\mathcal{X}_{\underline{k}}^t}$\,, which sends
$[z_{-1}:z_0:z_1]$ into :
\begin{equation}\label{eq:PROJ}
\Big[1:\frac{z_0}{z_{-1}}:\frac{z_1}{z_{-1}}:
\frac{z_1^{a_1}+tz_1^{k_1}}{z_{-1}^{a_1-1}z_0}:\cdots:
\frac{P_i}{z_{-1}^{w_i}z_0^{Z_{i-2}(a_2,\ldots,a_{i-1})}}:\cdots:
\frac{P_{r+1}}{z_{-1}^qz_0^{p-q}}\Big].
\end{equation}

Let $\rho'_{\underline{k}}:B'\PP^2\to \PP^2$ be the minimal
sequence of blow ups such that $\widehat{\pi}\circ
\rho'_{\underline{k}}$ extends to a regular map $B'\PP^2\to
\widehat{\mathcal{X}_{\underline{k}}^t}$. Let $L_\infty\subset
\PP^2$ be the line at infinity (defined by $z_{-1}=0$) and  $L_0$ the
closure in $\PP^2$ 
of $\{z_0=0\}$. We use the same notations for their strict
transforms via blow ups of $\PP^2$.

\begin{lemma}\label{prop:9.1}
$\widehat{\pi}\circ\rho'_{\underline{k}}$ sends $L_0$ and the
total transform of $L_\infty$ in $C^\infty_{\underline{k}}$.
\end{lemma}
\begin{proof} Use (\ref{eq:PROJ}) or the fact that the projection
  $\mbox{pr}:\mathcal{X}_{\underline{k}}^t
\to \C^2$ is regular and the corresponding restrictions of
$\mbox{pr}\circ(\widehat{\pi}\circ \rho'_{\underline{k}})$ and
$\rho'_{\underline{k}}$ are equal.
\end{proof}
Hence, from the point of view of $\mathcal{X}^t_{\underline{k}}$,
resolving the indeterminacy points of $\widehat {\pi}$ above
$L_\infty$ is irrelevant. Let $\rho_{\underline{k}}:B\PP^2\to
\PP^2$ be the {\it minimal} sequence of blow ups which resolve the
indeterminacies of $\widehat{\pi}$ {\it sitting in } $\C^2$ (hence
$\rho'_{\underline{k}}$ and $\rho_{\underline{k}}$ over $\C^2$
coincide). Denote by $E_\pi$ its exceptional curve and by $C_\pi$
the union of those irreducible components of $E_\pi$ which are
sent to $C^\infty_{\underline{k}}$.
% Set $B\C^2:=B\PP^2\setminus L_\infty$.
Summing up the above discussions, one obtains:

\begin{corollary}\label{cor:9.1} The restriction of \,
$\widehat{\pi}\circ\rho_{\underline{k}}$ induces an isomorphism
$B\PP^2\setminus (L_\infty\cup L_0\cup C_\pi)\to
\mathcal{X}_{\underline{k}}^t$. In particular, the Milnor fiber
can be realized as the complement of the projective curve
$L_\infty\cup L_0\cup C_\pi$ in $B\PP^2$.
\end{corollary}
\begin{proof} Use the fact that $\mathcal{X}_{\underline{k}}^t$ is smooth, cf.
(\ref{lem:9.2}).
\end{proof}

For the convenience of the reader, we represent in the following diagram
all the maps introduced in the previous discussion:

$$\xymatrix{
     B\PP^2  \ar[dr]_{\rho_{\underline{k}}}  
      & B'\PP^2 \ar[l] \ar[r]  \ar[d]^{ \rho'_{\underline{k}} }  
      &  \widehat{\mathcal{X}_{\underline{k}}^t} = 
            \mathcal{X}_{\underline{k}}^t \cup C^\infty_{\underline{k}}
      &   \mathcal{X}_{\underline{k}}^t \ar[d]^{pr} \ar@{_{(}->}[l] \\
    & \PP^2 = \C^{2} \cup L_{\infty} \ar@{-->}[ru]_{\hat{\pi}} 
    &
    &  \C^2 \ar@{_{(}->}[ll] }$$

\subsection{The curve-configurations $E_\pi$ and
$C_\pi$.}\label{ss:9.4}  The equations (\ref{eq:EGY}) and (\ref{eq:9.7}) show 
that the indeterminacy points of $\widehat{\pi}$ sitting in $\C^2$
are given by $\{z_0=P_{r+1}=0\}$. By equations (\ref{eq:9.6}) and
induction, this set equals
$\{z_0=P_2=0\}=\{z_0=z_1^{a_1}+tz_1^{k_1}=0\}$ sitting in $L_0$.
The indeterminacy at the points $(0,\xi_j)$, where $\{\xi_j\}_j$
are the roots of $z_1^{a_1-k_1}+t=0$, can be eliminated by a
single blow up (see e.g. below). The indeterminacy at $(0,0)$ (which
appears exactly when $k_1>0$, i.e. when $r>1$) requires, in
general, more blow-ups. The structure of $\widehat{\pi}$ at these
points will be revealed in the next paragraphs.

The modification $\rho_{\underline{k}}:B\PP^2\to \PP^2$ will be
constructed in two steps. First, we define a toric modification of
$\PP^2$ with exceptional curves $\cup_{j=2}^rV_j$\,, all above
$[1:0:0]$, such that $L_0\cup(\cup_{j=2}^rV_j)$ form a string.
After this modification, $\sum_{1 \leq i \leq r}(a_i -k_i)$
indeterminacy points survive, they will be eliminated in the
second step by blowing up each point once.

Recall that $z_0=\chi^{\overline{v}_0},
z_1=\chi^{\overline{v}_1}$. Denote by $(u_1,u_{r+1})\in N$ the
dual basis of $(\overline{v}_0, \overline{v}_1)$ and by
$\tilde{\sigma}$ the cone generated by it. Hence, the affine plane
$\C^2$ of coordinates $(z_0, z_1)$ is identified with the toric
surface $\mathcal{Z}_{\tilde{\sigma}, N}$. Take also
$u_0:=-(u_{r+1} + u_1)$ and the complete regular fan
$\mathcal{F}_0$ whose 1-dimensional cones are generated by $u_0,
u_1$ and $u_{r+1}$. Then $\mathcal{Z}_{\mathcal{F}_{0}, N}=\PP^2$.

\begin{figure}[ht!]
\labellist
\small\hair 2pt
\pinlabel $0$ at 228 159
\pinlabel $u_0$ at 137 72
\pinlabel $u_1$ at 435 265
\pinlabel $u_2$ at 352 288
\pinlabel $u_3$ at 350 396
\pinlabel $u_{r-1}$ at 154 396
\pinlabel $u_r$ at 78 275
\pinlabel $u_{r+1}$ at 98 208
\pinlabel $\sigma_0$ at 320 105
\pinlabel $\sigma_1$ at 432 342
\pinlabel $\sigma_2$ at 395 410
\pinlabel $\sigma_{r-1}$ at 46 408
\pinlabel $\sigma_r$ at 8 250
\pinlabel $\sigma_{r+1}$ at 28 106
\endlabellist
\centering
\includegraphics[scale=0.40]{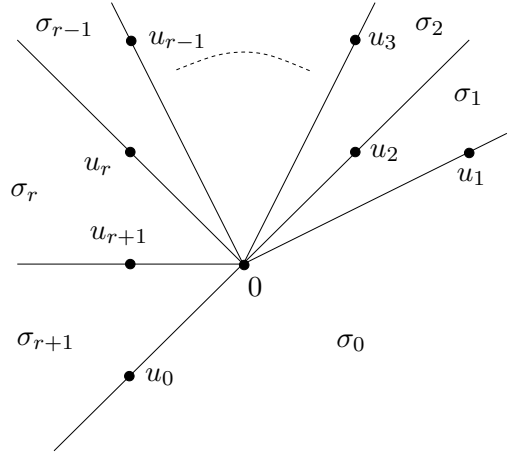}
\caption{The  complete regular fan $\mathcal{F}_{\underline{k}}$}
\label{fig:Fan}
\end{figure}

Next, consider the complete regular fan
$\mathcal{F}_{\underline{k}}$ subdividing $\mathcal{F}_0$,
whose 1-dimensional cones are generated by the primitive elements
$u_0,u_1,\ldots, u_{r+1}$ of $N$ such that (see Fig.
\ref{fig:Fan}):
\begin{equation} \label{fan}
   \left\{ \begin{array}{l}
               u_0+ u_2= (k_1 -1) u_1,\\
               u_{j-1} + u_{j+1} = k_j u_j \: \: \:   \mbox { for all } j \in
                \{2,\ldots,r\}, \\
               u_{r+1} + u_1 = - u_0.
           \end{array} \right.
\end{equation}
Its existence is ensured by the fact that $\underline{k}$ is an
admissible sequence which represents $0$.

$\mathcal{F}_{\underline{k}}$ being a subdivision of
$\mathcal{F}_0$, it induces a proper birational toric morphism:
\begin{equation} \label{propmorph}
   \mathcal{Z}_{\mathcal{F}_{\underline{k}}, N}
    \stackrel{\psi_{\underline{k}}}{\longrightarrow}
    \mathcal{Z}_{\mathcal{F}_{0}, N}=\PP^2.
\end{equation}
For any $j\in\{0,\ldots, r\}$, denote  by $O_j\simeq \C^*$ the
orbit in $ \mathcal{Z}_{\mathcal{F}_{\underline{k}}, N}$
corresponding to the 1-dimensional cone generated by $u_j$ and by
$V_j$ its closure. Clearly $V_0=L_\infty$ and $V_1=L_0$.

In fact,  $\psi_{\underline{k}}$ can also be
characterized independently of toric geometry: it is the unique
modification with exceptional divisors $\{V_j\}_{j=2}^r$, all
above the point $(0,0)$, such that $L_0\cup V_2\cup\cdots \cup
V_r$ form a string with self-intersections $1-k_1,\, -k_2,\ldots,
-k_r$, cf. \cite[(1.6)]{Oda}.

For all $j \in \{0,\ldots,r\}$ let $\sigma_j$ be the cone
generated by $(u_j, u_{j+1})$ and $(\xi_j, \eta_j)$ the dual basis
of $(u_j, u_{j+1})$. Set the corresponding monomials $x_j :=
\chi^{\xi_j}$ and $y_j:=\chi^{\eta_j}$. Therefore, in
$\mathcal{Z}_{\sigma_j, N}$ one has $\{x_j=0\} = V_j $ and
$\{y_j=0\} = V_{j+1}$.   
Moreover, the restriction:
\begin{equation} \label{partchart}
  \psi_{j}: \mathcal{Z}_{\sigma_j, N} \rightarrow
   \mathcal{Z}_{\tilde{\sigma}, N}=\C^2,
\end{equation}
of $\psi_{\underline{k}}$ to $\mathcal{Z}_{\sigma_j, N} \subset
\mathcal{Z}_{\mathcal{F}_{\underline{k}}, N}$
 is described by the following monomial changes of variables:
\begin{equation} \label{monom}
   \left\{ \begin{array}{l}
               z_{0}=x_j^{Z_{j-1}(k_1,\ldots,k_{j-1})}
                     y_j^{Z_j(k_1,\ldots,k_j)}\\
               z_{1} =
               x_j^{Z_{j-2}(k_2,\ldots,k_{j-1})}y_j^{Z_{j-1}(k_2,\ldots,k_j)}
           \end{array} \right.
   \ \mbox{ for all } j \in \{1,\ldots,r\}.
\end{equation}
This is a consequence of the fact that
for all $j \in \{1,\ldots,r\}$, one has:
\begin{equation}\label{eq:UUU}
u_j   = Z_{j-1}(k_1,\ldots, k_{j-1}) u_1 + Z_{j-2}(k_2,\ldots,
k_{j-1}) u_{r+1}.
  \end{equation}
(\ref{eq:UUU}) follows by (increasing) induction, (\ref{fan}), and
determinantal relations as in (\ref{recpol}).

\begin{theorem}\label{indetfin}
  Consider the birational map
  $\Psi_{\underline{k}}:=\psi_{\underline{k}}\circ \widehat{\pi}:
  \mathcal{Z}_{\mathcal{F}_{\underline{k}}, N} \dashrightarrow
  \widehat{\mathcal{X}_{\underline{k}}^t}$. The
  indeterminacy points of $\Psi_{\underline{k}}$ are contained in
  $\cup_{j=1}^r  O_j$. Moreover, each
  orbit $O_j$ contains precisely $a_j-k_j$ indeterminacy points, which are
  eliminated by blowing up each point once. Let $\rho_{\underline{k}}$ be
  the composition of
  $\Psi_{\underline{k}}$ with these blow ups, and $C$ any of the new
  exceptional curves
  obtained by one of these $\sum_{j=1}^r (a_j-k_j)$ blow ups. Then
  $\rho_{\underline{k}}(C)$
  is a curve with $\rho_{\underline{k}}(C)\not\subset
  C^\infty_{\underline{k}}$.

  Moreover,  $\Psi_{\underline{k}}$ (resp. $\rho_{\underline{k}}$)  maps \
  $\cup_{j=0}^r V_j$ (resp. their strict transforms) to
  $C^\infty_{\underline{k}}$.
\end{theorem}

\begin{proof}
Notice that for any $j\in\{1,\ldots,r\}$, in the chart
$\mathcal{Z}_{\sigma_j, N}$ with affine coordinates $(x_j,y_j)$
one has $\{x_j=0\}=O_j\cup\{V_j\cap V_{j+1}\}$, hence these affine
coordinate axes cover all the exceptional locus and indeterminacy
points. Hence, it is enough to analyze in each chart $(x_j,y_j)$
the behaviour of $\Psi_{\underline{k}}$ along $\{x_j=0\}$. Define
for each $i\in\{1,\ldots, r+1\}$ and $j\in\{1,\ldots, r\}$:
\begin{equation}\label{eq:MJI}
m_i^{(j)}:=\left\{\begin{array}{ll}
%Z_{i-2}(a_2,\ldots,a_{i-1})\cdot Z_{j-1}(k_1,\ldots,k_{j-1})+
Z_{j-i-1}(k_{i+1},\ldots,k_{j-1})& \mbox{if $i\leq j$}, \\ \ \\
%Z_{i-2}(a_2,\ldots,a_{i-1})\cdot Z_{j-1}(k_1,\ldots,k_{j-1})
-Z_{i-j-1}(a_{j+1},\ldots,a_{i-1})& \mbox{if $i> j$}.
\end{array}\right.
\end{equation}
The next technical lemma will not only guarantee that $m^{(j)}_i$
is the valuative order of $z_i\circ \psi_{\underline{k}}$ along
$V_j$, but it gives a rather complete structure of the pull-back
$z_i\circ\psi_j$ as well, where the maps $\psi_j$ are defined by
(\ref{partchart}):

\begin{lemma}\label{tech.lem} For any fixed $j$, one has:
\begin{equation}\label{eq:PPSI}
z_i\circ \psi_j=x_j^{m_i^{(j)}} \, y_j^{m_i^{(j+1)}}\, Q_i^{(j)}
\end{equation}
for some $Q_i^{(j)}\in\Z[t,x_j,y_j]$, which has  the following
properties too:

(a) $$Q_i^{(j)}\Big|_{x_j=0}=\left\{ \begin{array}{ll} c_1 &
\mbox{for $i\leq j$},\\ c_1' (c_2
y_j^{a_j-k_j}+c_3t)^{Z_{i-j-1}(a_{j+1},\ldots, a_{i-1})}&
\mbox{for $i>j$}\end{array}\right.$$ for some non-zero constants
$c_1, \ c_1', \ c_2$ and $c_3$, where $c_2$ and $c_3$ are
independent of $i$.

(b) Let $y_j=\xi$ be one of the roots of $c_2
y_j^{a_j-k_j}+c_3t=0$. For each $\xi$, expand  $Q_i^{(j)}$ in
Taylor series in local variables $(x_j,y_j-\xi)$, and write it as
a sum $\sum_{h\geq h_\xi} Q_i^{(j)}(h)$ of homogeneous polynomials
$Q_i^{(j)}(h)$ of degree $h$ in these local variables, such that
$Q_i^{(j)}(h_\xi)\not=0$. Then:
$$h_\xi=Z_{i-j-1}(a_{j+1},\ldots, a_{i-1}).$$
Hence, by (a) and (b), $x_j$ does not divide $Q_i^{(j)}(h_\xi)$.
\end{lemma}

\begin{proof}
The proof is straightforward and elementary. It uses for any fixed
$j\in \{1,\ldots, r\}$ induction over $i\in\{1,\ldots, r+1\}$, the
`inductive equations' (\ref{eq:9.2}), the substitution
(\ref{monom}), and inductive formulas relating $Z(\underline{x})$,
cf. (\ref{recpol}). For $i=1$, $z_1\circ\psi_j$ is given by
(\ref{monom}), which proves (\ref{tech.lem})  with $Q_1^{(j)}=1$.
The inductive step is given by (\ref{eq:9.2}), namely
\begin{equation}\label{eq:IND}\big(z_{i+1}\circ \psi_j\big)\cdot
 x_j^{m_{i-1}^{(j)}}y_j^{m_{i-1}^{(j+1)}}Q_{i-1}^{(j)}=
 \big(x_j^{m_i^{(j)}}y_j^{m_i^{(j+1)}}Q_i^{(j)}\big)^{a_i} +t
\big(x_j^{m_i^{(j)}}y_j^{m_i^{(j+1)}}Q_i^{(j)}\big)^{k_i}.
\end{equation}

 \noindent $\bullet$ {\em The case $1\leq i\leq j$.}  Assume
that for some $i\leq j-1$, (\ref{tech.lem}) is satisfied for both
$i$ and $i-1$. We will verify it for $i+1$. First we analyze in
(\ref{eq:IND}) the exponents of $x_j$ (the discussion for
$y_j$-exponents is similar). From the right hand side of
(\ref{eq:MJI}) one factors out $k_im_i^{(j)}$, and the inductive
step for these exponents which we need to verify is
$m_{i+1}^{(j)}=k_im_i^{(j)}-m^{(j)}_{i-1}$, which follows from
(\ref{recpol}). Next, the inductive formula for $Q^{(j)}_{i+1}$
is:
$$Q^{(j)}_{i+1} \cdot Q_{i-1}^{(j)}=
x_j^{(a_i-k_i)Z_{j-i-1}(k_{i+1},\ldots, k_{j-1})}
y_j^{(a_i-k_i)Z_{j-i}(k_{i+1},\ldots, k_{j})} (Q_i^{(j)})^{a_i}
+t(Q_i^{(j)})^{k_i}.$$ If $a_i>k_i$ then the exponent of $x_j$ is
positive, hence $Q_{i+1}^{(j)}\big|_{x_j=0}=
t(Q_i^{(j)}\big|_{x_j=0})^{k_i}\cdot
(Q_{i-1}^{(j)}\big|_{x_j=0})^{-1}$ is constant by induction. If
$a_i=k_i$ then one has a similar expression.

\vspace{1mm}

 \noindent $\bullet$ {\em The case $i=j+1$.} This is the first case when the
 $m_i^{(j)}$-expression 
changes its shape (cf. (\ref{eq:MJI})) and $Q_i^{(j)}\big|_{x_j=0}$
is not constant. Notice that  $m_j^{(j)}=0$ and $m_j^{(j+1)}=1$,
hence the inductive steps for the coordinate exponents can easily
be verified. Moreover, $Q_{j+1}^{(j)}\cdot
Q_{j-1}^{(j)}=y_j^{a_j-k_j}+t(Q_{j}^{(j)})^{k_j}$, hence
(\ref{tech.lem})(a-b) follows too with $h_\xi=1$.

\vspace{1mm}

 \noindent $\bullet$ {\em The case $i>j+1$.} The exponents of $x_j$ and $y_j$
 can be analyzed 
similarly, while  $Q^{(j)}_{i+1} \cdot Q_{i-1}^{(j)}=
(Q_i^{(j)})^{a_i} +t(Q_i^{(j)})^{k_i}\cdot M,$ where
$$M:=x_j^{(a_i-k_i)Z_{i-j-1}(a_{j+1},\ldots, a_{i-1})}
y_j^{(a_i-k_i)Z_{i-j-2}(a_{j+2},\ldots, a_{i-1})}.$$ Notice that
$Z_{i-j-1}(a_{j+1},\ldots,a_{i+1})$ is always strictly  positive (cf.
(\ref{lem:CC'})(3)). Hence, $M\big|_{x_j=0}=0$ if $a_i>k_i$ and
$=1$ otherwise. Hence (\ref{tech.lem})(a-b) follows again by
(\ref{recpol}).
\end{proof}
\noindent
The function  $z_i\circ \psi_j$ for $1\leq i\leq j$ is  regular,
while for $i>j$ it is
$$z_i\circ \psi_j=\frac{Q_i^{(j)}}{x_j^{Z_{i-j-1}(a_{j+1},\ldots,
    a_{i-1})}\, y_j^{
Z_{i-j-2}(a_{j+2},\ldots, a_{i-1})}}.$$ Notice that the exponent
$Z_{i-j-1}(a_{j+1},\ldots, a_{i-1})$ is always strictly  positive,
cf. (\ref{lem:CC'})(3). The 
$y_j$-coordinates of the indeterminacy points on $\{x_j=0\}$ are
given by $Q_i^{(j)}|_{x_j=0}$, which corresponds to the roots
$\xi$ introduced in the above technical lemma (\ref{tech.lem}). In
particular, by this lemma, any of them is eliminated by one blow
up.

All the other statements of Theorem (\ref{indetfin}) now follow easily.
This ends its proof.
\end{proof}

The previous theorem shows that $E_{\pi}$ has $(r-1) +
\sum_{i=1}^r(a_i -k_i)$ irreducible components and that $C_{\pi}=
\cup_{j=2}^r V_j$.

\begin{figure}[ht!]
\labellist
\small\hair 2pt
\pinlabel  $1$ at 5 120
\pinlabel  {$1-a_1$} at 112 120
\pinlabel $-a_2$ at 220 120
\pinlabel $-a_r$ at 473 120
\pinlabel $-1$ at 85 35
\pinlabel {$-1$} at 140 35
\pinlabel {$-1$} at 194 35
\pinlabel $-1$ at 248 35
\pinlabel $-1$ at 446 35
\pinlabel $-1$ at 500 35
\pinlabel $L_{\infty}=V_0$ at 5 90
\pinlabel $L_0=V_1$ at 160 90
\pinlabel $V_2$ at 246 90
\pinlabel $V_r$ at 500 90
\pinlabel {$a_1 -k_1$} at 114 -10
\pinlabel {$a_2 -k_2$} at 224 -10
\pinlabel {$a_r -k_r$} at 474 -10
\pinlabel $\dots$ at 114 50
\pinlabel $\dots$ at 224 50
\pinlabel $\dots$ at 474 50
\endlabellist
\centering
\includegraphics[scale=0.65]{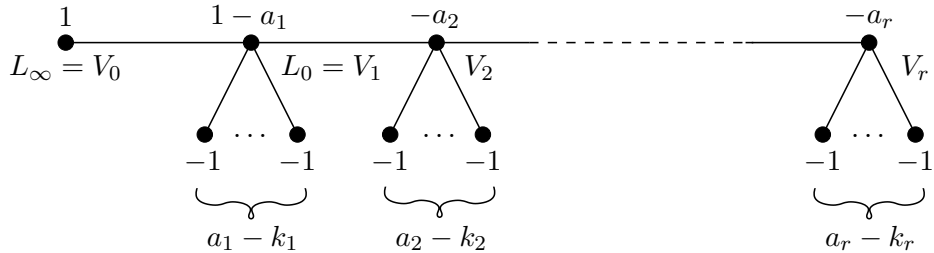}
\vspace{2mm}
\caption{Illustration for the Corollary (\ref{cor:9.2})}
\label{fig:Finalblow}
\end{figure}

Using the correspondence between the equations relating the
$u_i$'s in (\ref{fan}) and the self-intersections of the
corresponding curves in the associated toric variety
$\mathcal{Z}_{\mathcal{F}_{\underline{k}}, N}$,   we get:

\begin{corollary}\label{cor:9.2}
Consider the lines $L_\infty$ and $L_0$ on $\PP^2$ as above. Blow
up $r-1+\sum_{i=1}^r(a_i-k_i)$  infinitely close points of $L_0$
in order to get the  dual graph in Figure \ref{fig:Finalblow} of
the configuration of the total transform of  $L_\infty\cup L_0$
(this procedure topologically is unique, and its existence is
guaranteed by the fact that $\underline{k}\in K_r(\underline{a})$). Denote the
space obtained by this modification by $B\PP^2$. Then the Milnor
fiber $\mathcal{X}_{\underline{k}}^t$ of $S^{CS}_{\underline{k}}$
is diffeomorphic to $B\PP^2\setminus (\cup_{j=0}^r V_j)$.

Moreover, let $T$ be a small open tubular neighbourhood of
$\cup_{j=0}^r V_j$, and set
$F_{p,q}(\underline{k})=B\PP^2\setminus T$. Then
$F_{p,q}(\underline{k})$ is a representative of the Milnor fiber
of $S^{CS}_{\underline{k}}$ as a  manifold with boundary whose
boundary is $L(p,q)$.

Furthermore, the marking $\{V_i\}_i$ as in the Figure
\ref{fig:Finalblow}, defines on the boundary of
$F_{p,q}(\underline{k})$ an order; denote this supplemented space
by $F_{p,q}(\underline{k})^*$. Then its ordered boundary is
$L(p,q)^*$ endowed with the preferred order.
\end{corollary}

\subsection{The identification of Lisca's fillings with Milnor
fibers.}\label{ss:9.5} Let $W_{p,q}(\underline{k})^*$ be Lisca's
filling endowed with the preferred order on its boundary,  cf.
(\ref{ss:lisord});   and let $F_{p,q}(\underline{k})^*$ be the
Milnor fiber as in (\ref{cor:9.2}).

\begin{theorem}\label{th:9.1}
$W_{p,q}(\underline{k})^*$ is orientation-preserving diffeomorphic
to $F_{p,q}(\underline{k})^*$ by a diffeomorphism which preserves
the orders of the boundaries.
\end{theorem}
\begin{proof}
$F_{p,q}(\underline{k})^*$ from (\ref{cor:9.2}) satisfies the
criterion (\ref{glueplumb*}). Indeed, $B\PP^2$ with $V_0$
anti-blown-down differentiably, will serve as the differentiable
closed 4-manifold $V$. The homology  classes of the spheres $V_i$
are the $s_i$ ($ 1\leq i\leq r$), and the wished homology classes
$e$ with $e^2=-1$ are the classes of the $(-1)$ exceptional curves
from Figure \ref{fig:Finalblow}, multiplied by $\pm 1$. Moreover,
using the intersection form on $H_2(B\PP^2)$, we see that these
are the only classes $e$ with $e^2=-1$ which intersect
non-trivially only one component among $V_1,...,V_r$. In fact, all
the homological computations in $H_2(B\PP^2)$ fit perfectly with
Lisca's computation from \cite[\S 4]{L 08}. The compatibility of
orders is guaranteed by the compatibilities of the constructions,
see also (\ref{ss:4.3})-(\ref{markcoord}).
\end{proof}

\subsection{Remarks.}\label{finalrem} (1) Let $\rho_{\underline{k}}$ be the
modification introduced above (cf. (\ref{indetfin}) or
(\ref{ss:9.3})) (as the minimal modification which eliminates the
indeterminacy of $\widehat{\pi}|_{\C^2}$). Analyzing the proof of
(\ref{tech.lem}) we realize that $\rho_{\underline{k}}$ serves
also as the minimal modification which eliminates the
indeterminacy of the last component of $\pi$ from (\ref{eq:9.7}),
namely of the rational function $z_{r+1}= P_{r+1}/z_0^{p-q}.$ In
particular, we find the following alternative description of the
Milnor fiber $F_{p,q}(\underline{k})$:

\vspace{2mm}

 \emph{For each  $\underline{k} \in K_r(\underline{a})$, define the
  polynomial $P_{r+1}$ via the inductive system (\ref{eq:9.6}). Let
  $\rho_{\underline{k}}: B\PP^2\rightarrow \PP^2$ be the minimal
  modification of $\PP^2$ which eliminates the indeterminacy points of
$P_{r+1}/z_0^{p-q}$ sitting in $\C^2$. Then the dual graph of the
total transform of $L_{\infty} \cup L_0$ has the form indicated in
Figure \ref{fig:Finalblow}, and $F_{p,q}(\underline{k})$ is
orientation-preserving diffeomorphic to $B\PP^2 \setminus
(\cup_{j=0}^r V_j)$.}

 \vspace{2mm}

(2) One proves that the irreducible  decomposition of $P_{r+1}$
has the following form:
\begin{equation} \label{fineq}
               P_{r+1} =
     \prod_{j=1}^{r} \, \prod_{\ell=1}^{a_j-k_j}\, [ P_j +
     \xi_{j,\ell}\cdot z_0^{Z_{j-2}(a_2,\ldots,a_{j-1})}]
     ^{Z_{j-1}(k_{1},\ldots,k_{j-1})},
\end{equation}
where $\prod_{\ell=1}^{a_j-k_j} (\lambda +
     \xi_{j,\ell})=\lambda^{a_j-k_j}+t$.
     Moreover, the strict transforms
     by $\rho_{\underline{k}}$ of these irreducible components
     define `curvettas' of the $-1$ curves from Figure \ref{fig:Finalblow}.

\vspace{2mm}

(3) Sections (\ref{ss:9.1})--(\ref{ss:9.3}) contain some common
results with Balke's preprint \cite{B 99}. In fact, \cite{B 99}
convinced us that the identification (\ref{th:9.1}) should be
guided by a rather straightforward construction.

\medskip
\section{Invariants of 4-manifolds by closing
boundaries}
\label{capping}

In this section we present a procedure which provides invariants
for  $4$-manifolds with boundary $W$, by `closing' them with another
(fixed) 4-manifold $U$ (a `\emph{cap}'). Our main motivation is
Lisca's criterion
(\ref{glueplumb}) and his construction in \cite[\S 7]{L 08}. A
similar `closing'  will appear naturally for the Milnor fibers of
sandwiched singularities as well (see \ref{ss:11.1}).

Then we  generalize the results of section (\ref{ss:5.3}): we will
not only replace the plumbing 4-manifold $\Pi(\underline{a})$ by
an arbitrary 4-manifold $U$ (with the same boundary), but we will
show that the same criterion works for any such $U$ which
satisfies some  homological properties.

\subsection{The closing procedure.}\label{ss:10.1}
Let us  fix a 4-manifold with boundary $U$, which will be used as a
`cap' for other 4-manifolds $W$.

Assume that for some 4-manifold with boundary $W$, we have an
(orientation-preserving) diffeomorphism $\phi:\partial
W\rightarrow\partial \overline{U}$.  Then we construct the closed
manifold  $V=V(W,U,\phi)$  by gluing $W$ and $U$ along
their boundaries using $\phi$.    %  (see Figure \ref{fig:Cap}).
We say that $V$ is
obtained by \emph{closing} the boundary of $W$ by $U$. Its
diffeomorphism type  depends only on the isotopy class of $\phi$.
We write
$\mu: U \hookrightarrow V$  for the inclusion, and
$\partial_{U,\phi}$ for the composition $\phi^{-1}_*\circ
\partial_U:H_2(U,\partial U)\to H_1(\partial U)\to H_1(\partial
W)$.

Our goal is to establish some properties of $W$ read from the
homology of $V$.

In the sequel, we will suppose that \emph{$H_1(U)=0$ and $H_2(U)$ is
  free with a fixed  basis
$\underline{c}:= (c_1,\ldots,c_r)$}. Denote by
$\underline{c}^*:= (c_1^*,\ldots,c_r^*)$ the dual basis of $H_2(U,
\partial U)$ (via the intersection pairing
$H_2(U)\otimes H_2(U,\partial U)\to \Z$ which is a perfect pairing
under the above assumptions). Let $\mathcal{M}(\underline{c}):=
(Q_U(c_i, c_j))_{i,j}\in \mbox{Mat}_{r,r}(\Z)$ be the intersection matrix
of $U$.

Once we close $W$ by $U$, we concentrate on the following homological objects:
$\partial_{U,\phi}(\underline{c}^*):=
(\partial_{U,\phi}(c_1^*),\ldots,
\partial_{U,\phi}(c_r^*))$  in $H_1(\partial
W)^r$,  and the image $\mu_*(\underline{c}) \in H_2(V)^r$ of
$\underline{c}$.

\begin{proposition} \label{caprop}
Suppose that $\partial W$ is a rational homology sphere. Then, up to
an isomorphism (of such triplets), $(H_2(V), Q_{V}; \mu_*(\underline{c}))$
depends only on the manifold $W$, on
$\partial_{U,\phi}(\underline{c}^*) \in H_1(\partial W)^r$
and on $\mathcal{M}(\underline{c})$, but not
on the choice of the particular oriented 4-manifold $U$ (with
$H_1(U)=0$ and $H_2(U)$ free) used for the closing.
\end{proposition}

\begin{proof}
  The cohomological Mayer-Vietoris exact sequence, the vanishing
   $H^1(\partial W)=0$, and
   Poincar{\'e}-Lefschetz duality provide the exact sequence:
   $$0\rightarrow H_2(V)\rightarrow H_2(W,\partial W)\oplus H_2(U, \partial U)
   \stackrel{\Delta}{\rightarrow}H_1(\partial W).$$
    Hence  $H_2(V)=\ker  \Delta$, where  $\Delta(x\oplus y)=
     \partial_W(x)-\partial_{U,\phi}(y)$.  Consider the exact sequence:
  \begin{equation}\label{eq:10.1}
 0 \rightarrow H_2(U)\stackrel{i}{\longrightarrow}
 H_2(U,\partial U)\stackrel{\partial_{U,\phi}}
{\longrightarrow}
  H_1(\partial W).\end{equation}
The form $Q_U$ (given by  $\mathcal{M}(\underline{c})$) extends to
a rational form $Q_{U,\Q}$ on $H_2(U)_\Q$, and identifies
$H_2(U,\partial U)$ with the sublattice of elements $x\in
H_2(U)_\Q$  satisfying $Q_{U,\Q}(x,y)\in\Z$ for all $y\in H_2(U)$.
Hence, the restriction of $Q_{U,\Q}$ provides a rational form
$Q_{U,\partial U}:H_2(U,\partial U)^{\otimes 2}\to \Q$. In this
way we recover $H_2(U,\partial U)$ with its form $Q_{U,\partial
U}$  and the dual base $\underline{c}^*$,  and the sublattice
$H_2(U)$ in it. These, and the fact that $H_2(U)$ injects by
$y\mapsto (0\oplus i(y))$ into $\ker\Delta$, show that $H_2(V)$
and $\mu_*(\underline{c})\in H_2(V)^r$ can be recovered from the
input data.

Let us consider now $W$ instead of $U$. The analogue of sequence
(\ref{eq:10.1}) and a similar discussion as above show that the
form $Q_W$ extends to a rational form $Q_{W,\partial W}:
 H_2(W,\partial W)^{\otimes 2}\to \Q$. The point is that the wished
$Q_V$ is exactly the restriction of $Q_{W,\partial W}\oplus Q_{U,\partial U}$
on $\ker\Delta$ (which automatically takes only integral values).
\end{proof}

\subsection{The dependence on $\phi$.}\label{ss:10.2}
The following proposition shows that in the presence of an order,
in (\ref{caprop}) the choice of the gluing diffeomorphism $\phi$ is
irrelevant.

\begin{proposition}\label{prop:10.2}  Assume that $U$ is a
4-manifold with boundary such that $\partial \overline{U}$ is 
identified with $L(p,q)^*$, $H_1(U)=0$, and $H_2(U)$ is free  with
a fixed base $\underline{c}$.

Let $W$ be a Stein filling of $L(p,q)^*$ (i.e. on the boundary of
$W$ one can identify the preferred order of the lens space), and
let $V$ be obtained from $W$ by closing its boundary with $U$
using a gluing map which preserves the orientations and the orders
of the boundaries. Then $(H_2(V), Q_V; \mu_*(\underline{c}))$
(constructed in (\ref{caprop})) is independent of the choice of
$\phi$.

Moreover, $(H_2(V), Q_V;\mu_*(-\underline{c})) \simeq (H_2(V),
Q_V;\mu_*(\underline{c}))$ too.
\end{proposition}

\begin{proof} The argument is similar as in (\ref{ss:5.3}). The ambiguity
regarding $\phi$ stays in the group Diff$^{+,o}(L(p,q))$. If a
gluing   $\phi$ is replaced by  $\varphi\circ \phi$, where
$\varphi\in$Diff$^{+,o}(L(p,q))$ induces on $H_1(L(p,q))$
multiplication by $-1$, then we can twist $W$ by a
self-diffeomorphism which induces on the boundary $\varphi$ (as in
\ref{ss:5.3}), or instead, we can just multiply the homology of
$W$ by $-1$. The last isomorphism can be realized via
multiplication by $-1$ of $H_2(V)$.
\end{proof}

\medskip
\section{`From de Jong and van Straten to Lisca'}
\label{ljs}

\subsection{Closing the boundary of the Milnor fiber.}\label{ss:11.1}
 We keep all the notations of \S \ref{topsmooth}. We
consider again a decorated germ $(C,l)$ with smooth components $C_i$
and a picture deformation  $(C_S, l_S)$.

\begin{figure}[h!]
\labellist
\small\hair 2pt
\pinlabel $0'$ at 247 175
\pinlabel {$\overline{B'}$} at 185 173
\pinlabel $\overline{C_i'}$ at 306 225
\pinlabel $\sigma_i$ at 510 150
\pinlabel $T_i$ at 421 60
\endlabellist
\centering
\includegraphics[scale=0.40]{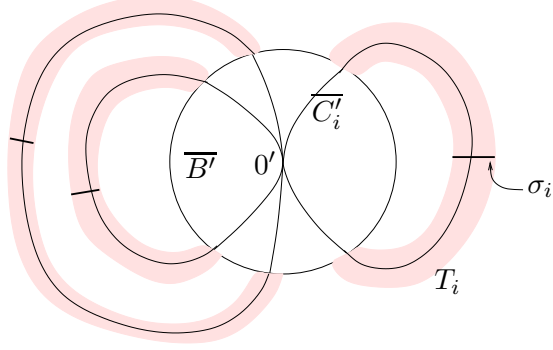}
\caption{The 4-manifold with boundary $U$}
\label{fig:Gluehandles}
\end{figure}

As the disc-configuration $D$ is obtained by deforming $C$,
 its  boundary $\partial D:= \linebreak \cup_{1\leq i \leq r}D_i
  \hookrightarrow \partial B$ is isotopic as an
  oriented link to $\partial C\hookrightarrow \partial B$. Therefore,
  we can isotope $D$
  outside a compact ball containing all the points $P_j$
  till its boundary  coincides with the boundary of $C$.
  Let $(B', C')$ be a second copy of $(B,C)$, and define:
  $$(V,\Sigma):=(B,D)\cup_{id}(\overline{B}', \overline{C}').$$
Here $V$ is
  the oriented 4-sphere obtained by gluing the boundaries
  of $B$ and $\overline{B}'$; and
  $\Sigma:= \cup_{i=1}^r\Sigma_i$, where
   $\Sigma_i$ is obtained by gluing
  $D_i$ (perturbed by the above isotopy) and $\overline{C}_i'$
  along their common boundaries.
Moreover,  one can also glue  $(\overline{B}', \overline{C}')$
with $(\tilde{B},\tilde{D})$ in such a way that the morphism
$\beta$ of  (\ref{blowpoints}) may be extended by the identity on
$\overline{B}'$, yielding:
$$ (\tilde{V}, \tilde{\Sigma})\stackrel{\beta}{\longrightarrow} (V,
\Sigma).$$ Here $\tilde{\Sigma}:= \cup_{i=1}^r \tilde{\Sigma}_i$,
where $\tilde{\Sigma}_i$ denotes the strict transform of the
sphere $\Sigma_i$, i.e. $\tilde{\Sigma}_i = \tilde{D}_i \cup
\overline{C}'_i$. Write $T:=\bigcup_{1\leq i \leq r} T_i$ and set
also  (see Figure \ref{fig:Gluehandles}):
\begin{equation} \label{defU}
  U:= \overline{B}'\cup T.
\end{equation}

Since $W=\tilde{B}\setminus T$ (cf. \ref{recofiber}), $\tilde{V}$
is obtained by closing the boundary of $W$ by the cap $U$. Our goal is to
recognize $W$ by a combination of  Lisca's criterion
(\ref{glueplumb}), of (\ref{caprop}), and of  (\ref{prop:10.2})
applied for this closing.  We start the needed preparations for
this program.

\begin{lemma} \label{handles} $U$ is independent of the chosen picture
deformation (therefore one may close all the different Milnor
fibers using the same $U$). In fact, each $T_i$ is a
$4$-dimensional handle of index $2$ glued to
    $\overline{B}'$ along the knot $\partial C_i\hookrightarrow
    \partial \overline{B}'$ endowed with the $(-l_i)$-framing.
  \end{lemma}

\begin{proof} As $H_2(V)=0$ (because $V\simeq \gras{S}^4$ ), we get
  $\Sigma_i^2 =0$. As $\tilde{\Sigma}_i$ is obtained from $\Sigma_i$
  by blowing-up
  (positively) $l_i$ points on it and taking its strict transform, we
  deduce that $\tilde{\Sigma}_i^2 =-l_i$. But this self-intersection
  is also equal to the self-linking number of the attaching circle
  $\partial \overline{C}_i'$ of the handle $T_i$ with respect to the
  attaching framing.
\end{proof}

Assume now that the decorated curve $(C,l)$ satisfying
$X(C,l)=\mathcal{X}_{p,q}$ is chosen as in \S \ref{cycl}. In
particular, $l$ is defined by (\ref{li}) or (\ref{la}).
We assume that the components of $C$ are marked as
in Figure \ref{fig:Diacurve}. We
write $W(\underline{a}, \underline{k})$ for $W$ and
$\tilde{V}(\underline{a}, \underline{k})=
 W(\underline{a}, \underline{k}) \cup
  U(\underline{a})$ for its closing by $U=U(\underline{a})$.

  \begin{lemma} \label{intsp}
     The intersection numbers of the oriented
     spheres $(\tilde{\Sigma}_i)_{1\leq i \leq r}$ inside the oriented
     $4$-manifold $\tilde{V}(\underline{a}, \underline{k})$
are the following:
     $$\left\{ \begin{array}{ll}
             \tilde{\Sigma}_i^2=-l_i & \mbox{for all} \
                    i \in \{1,\ldots,r\}, \\
             \tilde{\Sigma}_i \cdot \tilde{\Sigma}_j = 1-l_i &
             \mbox{for all}   \  i <j.
           \end{array} \right.$$
  \end{lemma}

\begin{proof}
  The first equalities were obtained for arbitrary decorated germs
  $(C,l)$ with smooth components $C_i$ during the proof of Lemma
  (\ref{handles}).
  For $i<j$, the surfaces $\tilde{\Sigma}_i$  and $\tilde{\Sigma}_j$
  meet at the origin of $\overline{B}'$. Therefore,
  $\tilde{\Sigma}_i\cdot_{\tilde{V}(\underline{a}, \underline{k})}
  \tilde{\Sigma}_j= \overline{C}_i\cdot_{\overline{B}} \overline{C}_j=
  - \overline{C}_i\cdot_{B} \overline{C}_j= -C_i\cdot C_j$.
  Then  apply  (\ref{inters}).
\end{proof}

  Consider the following homology classes in $H_2(U(\underline{a}))$:
  \begin{equation} \label{change}
     \left\{ \begin{array}{l}
          c_1:= [\tilde{\Sigma}_1]  \\
          c_i:= [\tilde{\Sigma}_i]-[\tilde{\Sigma}_{i-1}] \ \
          \mbox{for all} \
      \  i \in           \{2,\ldots,r\}.
     \end{array} \right.
   \end{equation}
A direct consequence of the previous lemma and of formula (\ref{la}) is:

  \begin{lemma} \label{bamboo}
     One has the following intersection numbers of the homology classes
     $c_i$:
     $$\left\{ \begin{array}{l}
             c_i^2=-a_i \ \ \mbox{for all} \  \ i \in \{1,\ldots,r\}, \\
             c_i \cdot c_j =
                 \left\{ \begin{array}{l}
                              1 \ \mbox{ if } \  |i-j|=1, \\
                              0 \ \mbox{ if } \  |i-j|>1.
                           \end{array} \right.
          \end{array} \right.$$
\end{lemma}

\noindent Next, we wish to identify $\partial_U(\underline{c}^*)$.

\begin{proposition} \label{eqfirst}
Using the notations (\ref{alphas}), one has a canonical
identification $\partial \overline{U}(\underline{a})\simeq
L(p,q)^*$ (i.e. which also identifies the preferred order on the
boundary).  Moreover (using the notations of (\ref{ss:4.new})), in
$H_1(\partial U)^r$  the following equality holds:
$$(\partial_U(c_1^*),\ldots, \partial_U(c_r^*))=\\
   \pm (\alpha^{\partial U}_1,\ldots, \alpha^{\partial U}_r).$$
\end{proposition}

\begin{proof}
 For each $i \in \{1,\ldots,r\}$ denote by
$\sigma_i$ a co-core of the handle $T_i$ (see Figure
\ref{fig:Gluehandles}) and orient it such that its intersection
number with $\tilde{\Sigma}_i$ is $+1$. Therefore
$([\sigma_1],\ldots,[\sigma_r])$ is the dual basis of
$([\tilde{\Sigma}_1],\ldots,[\tilde{\Sigma}_r])$ in
$H_2(U,\partial U)$. Hence, from equations (\ref{change}), we get:
\begin{equation} \label{derive}
   \left\{ \begin{array}{l}
          \partial_U[\sigma_i]=
          \partial_U(c_i^*)-\partial_U(c_{i+1}^*) \  \  \mbox{for all}
          \ \ i \in
                    \{1,\ldots,r-1\},\\

        \partial_U[\sigma_r]= \partial_U(c_r^*).
             \end{array}  \right.
\end{equation}

Look now at the oriented 4-manifold $\overline{U}$. By relation
(\ref{defU}), we see that $\overline{U}:= B'\cup (\bigcup_{1\leq i
\leq r} \overline{T}_i)$. We use the complex structure of $B'$ to
do blow-ups. Denote by
$\tilde{\overline{U}}\stackrel{\pi}{\rightarrow} \overline{U}$ the
composition of blow-ups of points above $0'\in \overline{U}$, such
that the dual graph of the preimage $\pi^{-1}(C')$ is isomorphic
to the one from Figure \ref{fig:Diacurve}. Then
$\tilde{\overline{U}}$ is a 4-manifold obtained by plumbing
according to the graph of Figure \ref{fig:Diarich} (this is
equivalent to (\ref{handles})).

\vspace{5mm}
\begin{figure}[ht!]
\labellist
\small\hair 2pt
\pinlabel $j_0$ at 5 130
\pinlabel  {$m_1 -1$} at 32 -20
\pinlabel $m_2$ at 210 -20
\pinlabel $m_t$ at 500 -20
\pinlabel $m_{t+1}$ at 681 -20
\pinlabel $\dots$ at 32 40
\pinlabel $\dots$ at 214 40
\pinlabel $\dots$ at 502 40
\pinlabel $\dots$ at 683 40
\pinlabel {$n_1 -3$} at 123 82
\pinlabel {$n_t -3$} at 593 82
\pinlabel $-2$ at 92 160
\pinlabel $-2$ at 151 160
\pinlabel $-2$ at 271 160
\pinlabel $-2$ at 439 160
\pinlabel $-2$ at 558 160
\pinlabel $-2$ at 618 160
\pinlabel $-1$ at 4 102
\pinlabel $-1$ at 63 102
\pinlabel $-1$ at 189 102
\pinlabel $-1$ at 247 102
\pinlabel $-1$ at 471 102
\pinlabel $-1$ at 533 102
\pinlabel $-1$ at 650 102
\pinlabel $-1$ at 715 102
\pinlabel $0$ at -5 70
\pinlabel $0$ at 74 70
\pinlabel $0$ at 174 70
\pinlabel $0$ at 256 70
\pinlabel $0$ at 459 70
\pinlabel $0$ at 545 70
\pinlabel $0$ at 634 70
\pinlabel $0$ at 726 70
\pinlabel {$-(m_1 +1)$} at 33 160
\pinlabel {$-(m_{2}+2)$} at 212 160
\pinlabel {$-(m_{t}+2)$} at 501 160
\pinlabel {$-(m_{t+1}+1)$} at 688 160
\pinlabel {$\tilde{\overline{\sigma}}_1$} at -10 32
\pinlabel {$\tilde{\overline{\sigma}}_{m_1-1}$} at 102 32
\pinlabel {$\tilde{\overline{\sigma}}_{m_1}$} at 162 32
\pinlabel {$\tilde{\overline{\sigma}}_{m_1 + m_2 -1}$} at 305 32
\endlabellist
*\vspace{5mm}
\centering
\includegraphics[scale=0.60]{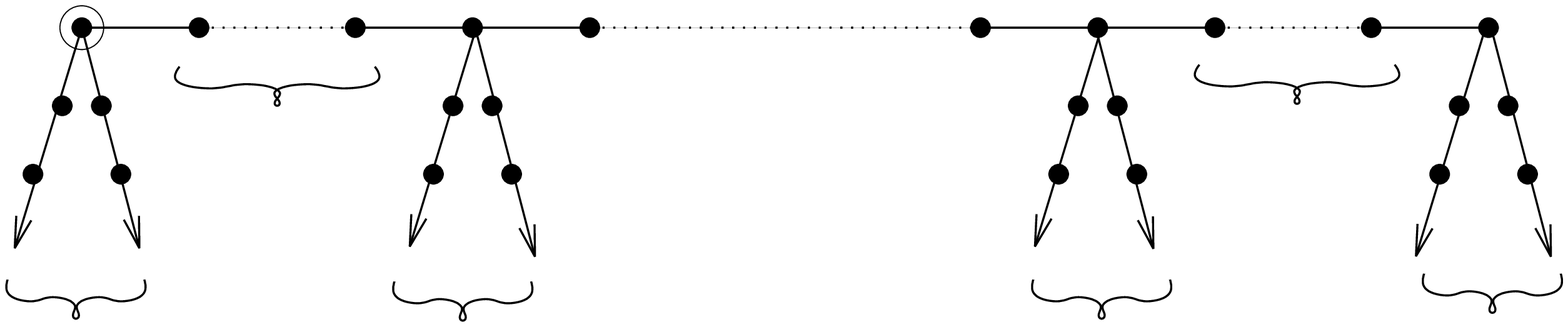}
 \caption{The plumbed 4-manifold $\tilde{\overline{U}}$}
\label{fig:Diarich}
\end{figure}

 Notice that its
boundary can canonically be identified with $L(p,q)$ (via plumbing
calculus). Indeed, first blowing-down the $(-1)$-curves and then
by anti-blowing-down (in the differential category) the
$(+1)$-curves arising from the $(0)$-curves, we get a plumbing
graph which without arrowheads is exactly the graph
$G(\underline{b})$. Considering in both graphs the preferred order
(cf. \ref{ss:4.3}), we get the proof of the first statement.
Notice that this also appoints the preferred order to $\partial
U(\underline{a})$.

In  Figure \ref{fig:Diarich}, the arrowheads denote again the
pre-images of the co-cores $\overline{\sigma}_i$ (i.e. $\sigma_i$
with opposite orientation;  they are analogs of the $R_i$'s of
(\ref{ss:4.3})). Therefore, by (\ref{repres}),  up to a
simultaneous change of sign, one has:
$$\partial_U[ \overline{\sigma}_i]=\nu_*(w_l) \ \  \mbox{ whenever } \
\  m_1 + \cdots +
m_{l-1} \leq i \leq m_1 + \cdots + m_l -1,$$ where the vectors
$w_l$ were defined in (\ref{precdual}).
Then from the  `duality relation' (\ref{precdual}):
$$ \partial_U [\sigma_i]= \alpha^{\partial U}_{i} - \alpha^{\partial
  U}_{i+1},$$
where $\alpha^{\partial U}_{r+1}:=\nu_*(\overline{v}_{r+1})=0$.
This combined with (\ref{derive}) ends the proof.
\end{proof}

If we sum up the results of this and the previous section, we get:
\begin{corollary}\label{cor:11}
Consider $\underline{s}\in H_2(\Pi(\underline{a}))^r$ as in
(\ref{glueplumb}), and $\underline{c}\in H_2(U(\underline{a}))^r$
defined in (\ref{change}). Then the following facts hold:

\vspace{1mm}

\noindent  \ \  I.\ (i) \
 $Q_{\Pi(\underline{a})}(s_i,s_j)=Q_{U(\underline{a})}(c_i,c_j)$
 for all $i,j$;

\vspace{1mm}

\ \ (ii) \ $\partial_{\Pi(\underline{a})}(\underline {s}^*)=\pm
\partial_{U(\underline{a})}(\underline {c}^*)$;

\vspace{1mm}

\noindent \ \ II. \   Let $W$ be a Stein filling of $L(p,q)^*$
(i.e. on the boundary of $W$ one can identify the preferred order
of the lens space), and close its boundary (using a
diffeomorphisms which preserves the orientations and the order of
the boundaries) by $\Pi(\underline{a})$ and $U(\underline{a})$
obtaining $V^\Pi$ resp. $V^U$. Then:
$$\big(H_2(V^\Pi),Q_{V^\Pi};\mu_*(\underline{s})\big)=
\big(H_2(V^U),Q_{V^U};\mu_*(\underline{c})\big).$$
\end{corollary}

This says that Lisca's criterion  (in order to recognize $W$),
expressed originally in $\big(H_2(V^\Pi),Q_{V^\Pi};\mu_*(\underline{s})\big)$
can be reinterpreted in
$\big(H_2(V^U),Q_{V^U};\mu_*(\underline{c})\big)$ too. Let us
apply this for the closing
$\tilde{V}(\underline{a}, \underline{k})=
 W(\underline{a}, \underline{k}) \cup
  U(\underline{a})$, and
search for the corresponding $(-1)$ curves. Set:
\begin{equation} \label{defE}
   E_j:=
  \beta^{-1}(P_j) \ \ \mbox{ for all }\  j \in \{1,\ldots,n\},
\end{equation}
where the number $n$ is defined in Theorem (\ref{expict}).

  \begin{proposition}\label{prop:11.11}
  One has the following equalities of matrices:
    $$ \begin{array}{l}
    \big(\tilde{\Sigma}_i\cdot E_j\big)_{1\leq i \leq r, 1\leq j \leq n}=
     \smallint D(\underline{a}; \underline{k}), \\
      \big(c_i\cdot E_j\big)_{1\leq i \leq r, 1\leq j \leq n }=
     D(\underline{a}; \underline{k}),
       \end{array}$$
 where  the entries of the left-side matrices are intersection numbers
  in  $H_2(\tilde{V}(\underline{a}, \underline{k}))$.
In particular, for any fixed $i\in\{1,\ldots , r\}$:
  $$ \#\{j \in \{1,\ldots, n\} \ \   | \ \
    c_i \cdot E_j \neq 0 \ \mbox{ but } \ c_k \cdot E_j =0 \
               \mbox{ for all } \ k \neq i    \} = a_i -k_i. $$
  \end{proposition}

\begin{proof}
By (\ref{incid}) the matrix $(\tilde{\Sigma}_i \cdot E_j)_{i,j}$
is equal to the incidence matrix of the picture deformation
corresponding to $\underline{k} \in K_r(\underline{a})$. Theorem
(\ref{expict}) implies the first equality of matrices. The second
one follows  from the construction of $\smallint D(\underline{a};
\underline{k})$ and from definition (\ref{change}). For the last
statement we search for columns of $ D(\underline{a};
\underline{k})$ with only one non-zero entry. For fixed $i$ they
correspond exactly to the block $M_{r,a_i-k_i}(i) $ of
(\ref{bloc}).
\end{proof}

Finally we get the searched isomorphism between the Milnor fibers
of the cyclic quotient singularity $\mathcal{X}_{p,q}$ and the
Stein fillings of the standard contact structure on $L(p,q)$:

\begin{theorem}\label{th:11}  Let $W(\underline{a},
\underline{k})^*$ be the Milnor fiber $W(\underline{a},
\underline{k})$ whose boundary is endowed with the preferred order
induced by the graph from Figure \ref{fig:Diarich} (which agrees
with the order of $\partial  \overline{U(\underline{a}})$ via the
gluing $\tilde{V}(\underline{a}, \underline{k}):=W(\underline{a},
\underline{k})\cup U(\underline{a})$). Then there is an
orientation-preserving diffeomorphism which preserves the orders
of the boundaries:
$$W(\underline{a}, \underline{k})^*\simeq W_{p,q}(\underline{k})^*.$$
\end{theorem}

\begin{proof} The statement follows from the criteria
(\ref{glueplumb}) and  (\ref{glueplumb*}) combined with
(\ref{cor:11}) once we check:
$$ \#\{e \in H_2(\tilde{V}(\underline{a}, \underline{k})) \  | \ e^2=-1, \
    c_i \cdot e \neq 0 \ \mbox{ but } \ c_k \cdot e =0 \
               \mbox{for all} \ k \neq i    \} = 2(a_i -k_i). $$
For this, first notice that $\tilde{V}(\underline{a},
\underline{k})$ is obtained from ${\mathbb S}^4$ by $n$ blow ups,
hence $\{[E_i]\}_{i=1}^n$ forms a basis in its second homology
group $H_2$, and the intersection form is diagonal with all
entries $-1$. This shows the equality of the sets $\{e\in H_2\ :\
e^2=-1\}=\{\pm [E_1],\ldots,\pm [E_n]\}$. Then use
(\ref{prop:11.11}).
\end{proof}

\medskip
\section{Final conclusions}
\label{concl}

\numberwithin{equation}{section}

 Let us list in short the two most important consequences
of the previous sections.

\begin{corollary}\label{cor:11.1}
  Once the order of the links  (or equivalently, the order of the coordinates
in the two constructions) are choosen in a compatible way,
  Christophersen and Stevens on one side and de Jong and van Straten
  on the other side parametrize in the same way the components of the
  miniversal base space of $\mathcal{X}_{p,q}$ by the elements
  $\underline{k}$ of $K_r(\underline{a})$:
  $$S_{\underline{k}}^{CS}= S_{\underline{k}}^{JS}.$$
\end{corollary}

\begin{corollary}\label{cor:11.2}
All the Milnor fibers $F_{p,q}(\underline{k})^*$ associated with
different smoothing components and endowed with the preferred
order on their ordered boundaries are different:  their boundaries
$L(p,q)^*$ and $\underline{k}\in K_r(\underline{a})$ determine
uniquely all the Milnor fibers up to orientation-preserving
diffeomorphisms which preserve the order of the boundary.
\end{corollary}


\begin{thebibliography}{00}

\bibitem{A 88} Arndt, J. \textit{Verselle Deformationen zyklischer
    Quotientensingularit{\"a}ten.} Diss. Hamburg, 1988.

\bibitem{A 74} Artin, M. \textit{Algebraic construction of Brieskorn's
    resolutions.} J. Algebra {\bf 29} (1974), 330-348. 

\bibitem{B 99} Balke, L. \textit{Smoothings of cyclic quotient
    singularities from a topological point of view.}
  arXiv:math/9911070

\bibitem{BHPV 04} Barth, W.P., Hulek, K., Peters, C.A.M., Van de Ven,
  A. \textit{Compact complex surfaces.} Second enlarged edition,
  Springer, 2004.

\bibitem{BR 95} Behnke, K., Riemenschneider, O. \textit{Quotient surface
    singularities and their deformations.} In \textit{Singularity theory.}
    D. T. L{\^e}, K. Saito \& B. Teissier eds. World Scientific, 1995,
    1-54.

\bibitem{BO 97} Bogomolov, F.A., de Oliveira, B. \textit{Stein Small
    Deformations of Strictly Pseudoconvex Surfaces.} Contemporary
  Mathematics \textbf{207} (1997), 25-41.

\bibitem{B 83} Bonahon, F. \textit{Diff{\'e}otopies des espaces
    lenticulaires.} Topology \textbf{22} (1983), 305-314.

\bibitem{CP 04} Caubel, C., Popescu-Pampu, P. \textit{On the contact
    boundaries of normal surface singularities.}
  C. R. Acad. Sci. Paris, Ser. I \textbf{339} (2004) 43-48.

\bibitem{CNP 06} Caubel, C., N{\'e}methi, A., Popescu-Pampu,
  P. \textit{Milnor open books and Milnor fillable contact
    3-manifolds.}   Topology \textbf{45} (2006), 673-689.

\bibitem{C 91} Christophersen, J.A. \textit{On the components and discriminant
   of the versal base space of cyclic quotient singularities.} In
 \textit{Singularity
   theory and its applications.} Warwick 1989, Part I, D. Mond,
 J. Montaldi eds.,
   LNM \textbf{1462}, Springer, 1991.

\bibitem{E 90} Eliashberg, Y. \textit{Filling by holomorphic discs and
    its applications.}   Geometry of low-dimensional manifolds, 2
  (Durham, 1989),  45-67, London Math. Soc. Lecture Note Ser., \textbf{151},
  Cambridge Univ. Press, 1990.


\bibitem{GS 99} Gompf, R.E., Stipsicz, A.I. \textit{4-manifolds and
    Kirby calculus.} A.M.S., 1999.

\bibitem{G 72} Grauert, H. \textit{{\"U}ber
    die Deformationen isolierter Singularit{\"a}ten analytische Mengen.}
  Invent. Math. \textbf{15} (1972), 171-198.

\bibitem{JS 98} de Jong, T., van Straten, D. \textit{Deformation theory of
  sandwiched singularities.} Duke Math. Journal \textbf{95}, No. 3 (1998),
  451-522.

\bibitem{K 85} Koll{\'a}r, J. \textit{Toward moduli of singular varieties.}
   Compositio Math. \textbf{56} (1985), 369-398.

\bibitem{KS 88} Koll{\'a}r, J., Shepherd-Barron,
  N. I. \textit{Threefolds and deformations of surface singularities.}
  Invent. Math. \textbf{91} (1988), 299-338.

 \bibitem{L 00} L{\^e}, D.T. \textit{Les singularit{\'e}s sandwich.} In
   \textit{Resolution of singularities: a research textbook in tribute
   to Oscar Zariski.} H. Hauser, J. Lipman, F.Oort eds,
   Birkh{\"a}user, 2000,   457-483.


  \bibitem{L 04} Lisca, P. \textit{On lens spaces and their symplectic
      fillings.}
    Math. Res. Letters \textbf{1}, vol. 11 (2004), 13-22.
    arXiv:math.SG/0203006.

\bibitem{L 08} Lisca, P. \textit{On symplectic fillings of lens spaces.}
Trans. Amer. Math. Soc. \textbf{360} (2008), 765-799.
 arXiv:math.SG/0312354.

\bibitem{L 84} Looijenga, E. J. N. \textit{Isolated singular points on
    complete intersections.} London Mathematical Society Lecture Note
  Series, {\bf 77}. Cambridge University Press, 1984. 

\bibitem{M 90} McDuff, D. \textit{The structure of rational and ruled
    symplectic $4$-manifolds.}   J. Amer. Math. Soc.  \textbf{3}
  (1990),  no. 3, 679--712.

\bibitem{M 03} M{\"o}hring, K. \textit{On sandwiched singularities.}
   Thesis, Univ. Mainz, 2003.


\bibitem{N 81} Neumann, W. \textit{A calculus for plumbing applied to
    the topology of complex surface singularities and degenerating
    complex curves.} Trans. Amer. Math. Soc. \textbf{268}, 2 (1981), 299-344.

\bibitem{Oda} Oda, T. \textit{Convex bodies and algebraic
geometry}. Springer-Verlag, 1988.

\bibitem{OO 03} Ohta, H., Ono, K. \textit{Symplectic fillings of the
    link of simple elliptic singularities.} J. Reine
  Angew. Math. \textbf{565} (2003), 183-205.

\bibitem{OO 05} Ohta, H., Ono, K. \textit{Simple singularities and
    symplectic fillings.} J. Differential Geom. \textbf{69} (2005),
  1-42.

 \bibitem{OW 77} Orlik, P., Wagreich, P. \textit{Algebraic surfaces with
 $k^*$-action.} Acta Math. \textbf{138} (1977), 43-81.



\bibitem{PPP 07} Popescu-Pampu, P. \textit{The geometry of
continued fractions and the topology of surface singularities.}
Advanced Stud. in Pure Maths \textbf{46} (2007), 119-195.

\bibitem{R 74} Riemenschneider, O. \textit{Deformationen von
    Quotientensingularit{\"a}ten (nach Zyklischen Gruppen).}
    Math. Ann. \textbf{209}, 211-248 (1974).


\bibitem{S 68} Schlessinger, M. \textit{Functors of Artin rings.}
  Trans. Amer. Math. Soc. \textbf{130} (1968), 208-222

\bibitem{S 90} Spivakovsky, M. \textit{Sandwiched singularities and
   desingularization of surfaces by normalized Nash transformation.}
   Ann. Math. (2) \textbf{131} (1990), 411-491.


\bibitem{S 91} Stevens, J. \textit{On the versal deformation of cyclic
  quotient singularities.} In \textit{Singularity
   theory and its applications.} Warwick 1989, Part I, D. Mond,
 J. Montaldi eds., LNM \textbf{1462}, Springer, 1991, 312-319. 

\bibitem{S 03} Stevens, J. \textit{Deformations of singularities.}
  Springer LNM \textbf{1811}, 2003.

\bibitem{W 81} Wahl, J. \textit{Smoothings of normal surface singularities.}
   Topology \textbf{20} (1981), 219-246.

\end{thebibliography}
\end{document}